\documentclass[11pt]{elsarticle} 

\makeatletter
\def\ps@pprintTitle{%
	\let\@oddhead\@empty
	\let\@evenhead\@empty
	\def\@oddfoot{}%
	\let\@evenfoot\@oddfoot}
\usepackage{bm}
\usepackage{bbm}
\usepackage{amsmath}
\usepackage{amsthm}
\usepackage{empheq}
\usepackage{float}
\usepackage{indentfirst}
\usepackage{cases}
\usepackage{algorithm}  
\usepackage{algorithmicx} 
\usepackage{algpseudocode} 
\usepackage{subcaption}
\usepackage{graphicx}
\usepackage{booktabs}
\usepackage{diagbox}
\usepackage[all]{xy}
\usepackage{indentfirst}
\usepackage{epstopdf}
\usepackage{epsfig}
\usepackage{overpic}
\usepackage{array}
\usepackage{color}
\usepackage{multirow}
\usepackage[mathscr]{euscript}
\usepackage{latexsym,bm}
\usepackage{tikz}
\usetikzlibrary{shapes,arrows}
\usepackage[left=2.5cm,right=2.5cm,top=2cm,bottom=2cm]{geometry}
\usepackage{enumerate}
\usepackage{physics}

\usepackage{ulem}

\usepackage{framed,multirow}

\usepackage{amssymb}
\usepackage{latexsym}
\usepackage{mathtools}

\usepackage{url}
\usepackage{xcolor}
\definecolor{newcolor}{rgb}{.8,.349,.1}

\renewcommand{\vec}[1]{\boldsymbol{#1}} 

\newtheorem{assumption}{Assumption}

\def\bx{\bm{x}}

\def\by{\bm{y}}

\def\bb{\bm{b}}
\def\me{\mathbbm{e}}
\def\D{\textup{d}}
\def\mone{\mathbbm{1}}
\def\Var{\mathrm{Var}}

\makeatletter

\newcommand{\Rmnum}[1]{\expandafter@slowromancap\romannumeral #1@}
\makeatother

\tikzset{global scale/.style={
		scale=#1,
		every node/.append style={scale=#1}
	}
}
\tikzstyle{startstop} = [rectangle,rounded corners, minimum width=3cm,minimum height=1cm,text centered, draw=black,fill=red!30]
\tikzstyle{io} = [trapezium, trapezium left angle = 70,trapezium right angle=110,minimum width=3cm,minimum height=1cm,text centered,draw=black,fill=blue!30]
\tikzstyle{process} = [rectangle,minimum width=3cm,minimum height=1cm,text centered,text width =3cm,draw=black,fill=orange!30]
\tikzstyle{decision} = [diamond,minimum width=1.7cm,minimum height=0.5cm,shape aspect=2, text centered,draw=black,fill=green!30]
\tikzstyle{arrow} = [thick,->,>=stealth]

\newtheorem{definition}{Definition}
\newtheorem{theorem}{Theorem}
\newtheorem{lemma}[definition]{Lemma}
\newtheorem{corollary}[definition]{Corollary}
\begin{document}

\begin{frontmatter}

\title{\LARGE Stochastic particle method with birth-death dynamics}

\author[PKU]{Jingyang Huang}
\ead{jyhuang03@stu.pku.edu.cn}

\author[PKU]{Zhengyang Lei}
\ead{leizy@stu.pku.edu.cn}

\author[PKU]{Sihong Shao}
\ead{sihong@math.pku.edu.cn}


\address[PKU]{CAPT, LMAM and School of Mathematical Sciences,  Peking University, Beijing 100871, China}


\begin{abstract}
In order to numerically solve high-dimensional nonlinear PDEs and alleviate the curse of dimensionality, a stochastic particle method (SPM) has been proposed to capture the relevant feature of the solution through the adaptive evolution of particles [J. Comput. Phys. 527 (2025) 113818]. In this paper, we introduce an active birth-death dynamics of particles to improve the efficiency of SPM. The resulting method, dubbed SPM-birth-death, sample new particles according to the nonlinear term and execute the annihilation strategy when the number of particles exceeds a given threshold. A rigorous error estimation for SPM-birth-death is established, elucidating the first-order convergence in time and space, as well as half-order accuracy in the initial sample size with explicit variance estimates. We also extend the analysis framework to SPM and provide theoretical justification for the existing numerical convergence study. 
Our theoretical results reveal that the introduced active birth-death dynamics of particles results into less frequent resampling and SPM-birth-death is thus able to achieve higher efficiency than SPM. Validating benchmarks are provided. In particular, preliminary numerical experiments on the Allen-Cahn equation demonstrate that SPM-birth-death can achieve smaller errors at the same computational cost compared with the original SPM.
\end{abstract}

\begin{keyword}
Stochastic particle method; 
Birth-death process;
Markovian property;
Conditional independence;
Error estimation;
Adaptive sampling;
Piecewise constant reconstruction;  
Allen-Cahn equation
\end{keyword}

\end{frontmatter}

\section{Introduction}
With advances in computational mathematics, numerical methods for solving high-dimensional partial differential equations (PDEs) have attracted increasing attention from researchers. Typical high-dimensional PDEs include the 6-D Vlasov equation \cite{kormann2019massively, vlasov1961} and the Boltzmann equation \cite{dimarco2018efficient, bird1970direct, cercignani1988boltzmann}, which describe the evolution of charged particles in plasma physics, as well as the Hamilton–Jacobi–Bellman equation \cite{bk:Bellman1957, dolgov2021tensor} in control theory, whose dimension is determined by the number of players. The three well-known traditional numerical methods—finite difference, finite element, and spectral methods—are highly mesh-dependent. Consequently, they suffer from the curse of dimensionality (CoD) \cite{bk:Bellman1957} when applied to high-dimensional PDEs; that is, their computational cost increases exponentially with dimension.

Several numerical attempts have been made to alleviate the CoD, yet these approaches often suffer from limitations and drawbacks. When solutions of high-dimensional PDEs exhibit a low-rank structure, the tensor train method can exploit this structure via low-rank tensor grids, significantly reducing storage costs \cite{bachmayr2023low,bachmayr2016tensor,dolgov2021tensor, tang2024solving, richter2024continuous}. However, the effectiveness of this method hinges on the low-rank assumption, which is not always valid, thereby limiting its applicability. The sparse grid method \cite{Griebel2004Sparsegrids, Griebel1990ACT, Griebel2004Schrodinger} reduces storage requirements by applying a hyperbolic truncation to the tensor product basis set and discarding coefficients with absolute values below a given tolerance. While effective, this coefficient omission can compromise numerical stability \cite{kormann2016sparse}, and the accuracy is highly sensitive to solution smoothness. Deep-learning-based PDE solvers, currently undergoing rapid development, have been successfully applied to a wide range of high-dimensional problems \cite{HanJentzenE2018,RaissiPerdikarisKarniadakis2019,HurePhamWarin2020,Weinancontinuous,sirignano2018deep,e2018deepritz,li2021fourier,lu2021deeponet,Zeng2022adaptive}. However, such solvers exhibit high sensitivity to hyperparameter tuning \cite{kast2024positional}, posing challenges for theoretical understanding \cite{richter2021solving}.

Leveraging the Monte Carlo method's weak dimensionality dependence, particle methods have emerged as a promising direction for solving high-dimensional problems. Notable examples include: the Particle-in-Cell (PIC) method for the Vlasov equation \cite{chen2023fleks,luu2016voronoi,Teunissen2014}, branching random walk algorithms for the Wigner equation \cite{ShaoXiong2020, XiongShao2020Overcoming}, and branching diffusion processes for semi-linear parabolic equations \cite{HenryLabordere2019, chorin2009stochastic}. Recently, a stochastic particle method (SPM) was proposed for moderately high-dimensional nonlinear PDEs \cite{Lei2025}. This approach simulates solution evolution through particle motion and weight updates while maintaining a constant particle count throughout computation. This paper introduces an SPM incorporating a birth-death mechanism, dubbed SPM-birth-death,  that adaptively generates particles in regions requiring enhanced resolution to improve accuracy, and annihilates them in over-resolved areas to boost efficiency.

    Prior to this work, particle generation and annihilation mechanisms have been explored in the literature. The hybrid method with deviational particles (HDP) proposed in \cite{Yan2015, Yan2016} generates new particles from the source term at each time step and adds them to the original particle system. Consequently, HDP incorporates a particle annihilation mechanism to prevent exponential growth in particle number. This is implemented by first discarding the original particles, followed by resampling a reduced set from the numerical solution. In the Wigner branching random walk algorithm \cite{ShaoXiong2020, XiongShao2020Overcoming}, a single particle splits into three upon jumping, while adopting an annihilation strategy that cancels opposing-sign weighted particles within close proximity in phase space. In both methods, particle birth/death constitutes a passive, algorithmically essential component. Critically, the proposed SPM-birth-death actively incorporates a birth-death mechanism into particle methods to enhance computational efficiency. Within plasma simulation PIC methods \cite{chen2023fleks,luu2016voronoi,Teunissen2014}, active particle management strategies have been documented—for instance: splitting heavy-weight particles into light-weight counterparts to maintain accuracy, and merging light-weight particles into heavy-weight ones to reduce computational costs.


The rest of this paper is organized as follows. In Section 2, we provide a brief overview of SPM originally proposed in \cite{Lei2025}. Section 3 details the particle generation and annihilation mechanisms, leading to the construction of SPM-birth-death method. A rigorous error analysis for SPM-birth-death and SPM is presented in Section 4, followed by complexity analysis and validating benchmarks in Section 5. Finally, Section 6 concludes the paper with a discussion on future work.

\section{SPM: A brief}\label{sec:review_SPM}

Numerical resolution for high-dimensional nonlinear PDEs is a challenging task due to the curse of dimensionality. Recently, the stochastic particle method (SPM) proposed in \cite{Lei2025} provides an efficient way to solve general nonlinear PDEs of the form
\begin{equation}\label{cauchy problem}
	\begin{aligned}
		\frac{\partial } {\partial t}  u(\bx, t) &= \mathcal{L} u(\bx, t) +  f(t,\bx,u,\vec{\nabla }u), \quad \bx \in \mathbb{R}^d, ~ t> 0, \\
		u(\bx, 0) &= u_0(\bx), \quad \bx \in \mathbb{R}^d, 
	\end{aligned}
\end{equation}
where $\mathcal{L}$ represents a linear operator, $f(t,\bx,u,\vec{\nabla }u)$ is the nonlinear term, and the initial data $u_0(\bx) \in L^2(\mathbb{R}^d)$. Typical examples of $\mathcal{L}$ include the advection operator $\vec{b}\cdot \vec{\nabla }$ for transport equations, the Laplacian operator $c\Delta$ for diffusion equations and nonlocal operator $-(-\Delta)^{\alpha/2}$. SPM has been successfully applied to solve some typical high-dimensional nonlinear PDEs, such as the Allen-Cahn equation and Hamilton-Jacobi-Bellman equation, demonstrating its effectiveness and efficiency.

The core idea of SPM is to approximate the solution at the instant $t_m$ in the weak sense, using a set of $N$ particles with positions and real-valued weights. All particles at $t_m$ form a weighted point distribution
\begin{equation}\label{Xt_original}
    X_{t_m} := \frac{1}{N}\sum_{i=1}^{N} w_i(t_m) \delta_{\bx_i(t_m)},
\end{equation}
where $\bx_i(t_m)\in \mathbb{R}^d$ and $w_i(t_m)\in\mathbb{R}$ denote the position and weight of the $i$-th particle at time $t_m$, respectively, and $\delta_{\bx}$ is the Dirac delta function centered at $\bx$. The distribution approximates the continuous solution $U_m(\bx)\approx u(\bx, t_m)$ such that for any suitable test function $\varphi(\bx)$, the following holds:
\begin{equation}\label{eq:weak formulation}
		\langle \varphi, U_m\rangle
		\approx \langle \varphi,  X_{t_m} \rangle 
		= \frac{1}{N} \sum_{i=1}^{N} w_i(t_m) \varphi(\bx_i(t_m)),
\end{equation}
where $\langle f, g\rangle  = \int_{\mathbb{R}^d} f(\bx) g(\bx) \D \bx$ is the standard inner product, the set of test functions is chosen as $\varphi(\bx) = \mone_{Q_{k}}(\bx), k = 1,2,\dots,K$, $\mone_{Q_{k}}(\bx)$ denotes the indicator function of sub-domain $Q_k$ and $\Omega = \bigcup_{k=1}^K Q_k$ constitutes a partition of the computational domain $\Omega$. The weak formulation~\eqref{eq:weak formulation} defines a high-dimensional integral that depends on time and is calculated using Monte Carlo methods. Hence, the idea of  importance sampling can be used to reduce the variance, suggesting that the distribution of particle locations $p(\bx)$ at $t_m$ should be chosen to capture the main features of the solution $U_m(\bx)$, or more specifically, $p(\bx)=\frac{|U_m(\bx)|}{\int_{\mathbb{R}^d} |U_m(\bx)| \D \bx}$ with weight $w_i(t_m)=\frac{U_m(\bx_i(t_m))}{p(\bx_i(t_m))}$.

The central task of SPM is to evolve the particle system $\{\bx_i(t_m),w_i(t_m)\}$ to approximate the solution $U_{m}(\bx)$ at each time step, where existing time discretization methods naturally provide the guiding principle. The authors in \cite{Lei2025} take Lawson-Euler scheme \cite{lawson1967generalized} as an example to establish the evolution rules of $u(\bx,t)$ in the time direction,
\begin{equation}\label{eq:lawson}
	U_{m+1}(\bx) = \me^{\tau \mathcal{L}} \left(U_m(\bx) + \tau  f(t_m,\bx,U_m,\vec{\nabla }U_m)\right),
\end{equation}
where $\tau$ is the time step and $U_m(\bx)$ denotes the numerical solution. Thus the particle motion and weight update rules can be derived as 
\begin{subequations}\label{eq:SPM}
	\begin{empheq}[left={(\textbf{SPM})\quad\empheqlbrace}]{align}
		\label{eq:B_loc}&\bx_i(t_m)\xrightarrow{\text{relocationg}}\tilde{\bx}_i(t_m)\xrightarrow{\me^{\tau\mathcal{L}^*}} \bx_i(t_{m+1}),\\
		\label{eq:B_weight}&w_i(t_{m+1}):= \frac{[U_m(\tilde{\bx}_i(t_m)) + \tau f(t_m,\tilde{\bx}_i(t_m),U_m,\nabla U_m)]Z_m}{|U_m(\tilde{\bx}_i(t_m)) + \tau f(t_m,\tilde{\bx}_i(t_m),U_m,\nabla U_m)|}.
	\end{empheq}
\end{subequations}
where the relocation step is to redistribute particles $\{\tilde{\bx}_i(t_m)\}_{i=1}^N$ according to
\begin{equation}
	\label{eq:relocating}
	\tilde{\bx}_i(t_m) \sim \frac{|U_m(\bx) + \tau f(t_m,\bx,U_m,\nabla U_m)|}{Z_m}, \quad Z_m = \int_{\mathbb{R}^d} |U_m(\bx) + \tau f(t_m,\bx,U_m,\nabla U_m)| \D \bx, 
\end{equation}
and $\mathcal{L}^*$ denotes the adjoint operator of $\mathcal{L}$ defined on $\mathbb{H}^2(\mathbb{R}^d)$. That is, the particle first relocates to a new position $\tilde{\bx}_i(t_m)$ according to the updated solution magnitude, and then moves under the dynamics governed by the semigroup $\me^{\mathcal{L}^*\tau}$. In~\cite{Lei2025}, the authors present two strategies and recommend the aboved SPM because it is designed to respect the importance sampling principle.

At each time step, SPM requires the evaluation of $U_m(\bx)$ and $\nabla U_m(\bx)$ at particle locations to update weights, see Eq.~\eqref{eq:B_weight}. However, since the solution is represented by a weighted point distribution in Eq.~\eqref{Xt_original}, a reconstruction process is necessary to obtain continuous approximations of $U_m(\bx)$ and $\nabla U_m(\bx)$ from the discrete particle data. In \cite{Lei2025}, the so-called virtual uniform grid (VUG) is employed to approximate $U_m(\bx)$ by the piecewise constant function. Let $\Omega = \prod\limits_{j=1}^d \left[l_j,r_j\right] \subset \mathbb{R}^d$ denote the computational domain containing all particles. Suppose the particle locations are given by
	\begin{equation*}
		\bx_i(t_m) = \left(x_i^1,x_i^2,\dots,x_i^d\right), \quad i = 1,2,\dots,N(t_m),
	\end{equation*}
	and let $l_j = \min\limits_{1\leq i\leq N(t_m)} x_i^j$ and $r_j = \max\limits_{1\leq i\leq N(t_m)} x_i^j$. Then, the computational domain can be decomposed into several grids $Q_k$ with side length $h$,
	\begin{equation}\label{eq:decom}
		\Omega={\prod\limits_{j=1}^d \left[l_j,r_j\right]} = \bigcup_{k=1}^K Q_k.
	\end{equation}
	Suppose that the function value at location $\vec{x}$ needs to be estimated, and it lies within the grid $Q_k$. Then, based on the weak approximation~\eqref{eq:weak formulation} of the particle system, we have
	\begin{equation}\label{eq weight sum approx u}
		\begin{aligned}
			U_m(\bx) 
			&\approx  \mathcal{P}_h{U}_m(\bx) := \sum_{k=1}^K\left(\frac{1}{|Q_k|}\int_{Q_k} U_m(\by) \D \vec{y}\right) \mone_{Q_{k}}(\bx)\\
			&\approx \sum_{k=1}^K \left(\frac{1}{N|Q_k|}\sum_{i:x_i(t_m)\in Q_k}w_i(t_m)\right) \mone_{Q_{k}}(\bx)=\mathcal{P}_hX_{t_m}=:\tilde{U}_m.
		\end{aligned}
	\end{equation}
Here $\mathcal{P}_hU_m$ and $\mathcal{P}_hX_{t_m}$ denote the piecewise constant projection of $U_m$ and $X_{t_m}$, respectively. In fact, we can derive the second approximation in \eqref{eq weight sum approx u} by substituting the test function $\varphi(\bx)$ with the indicator function $\mone_{Q_{k}}(\bx)$ in \eqref{eq:weak formulation}.
	Once the estimation of the function value is obtained, the derivative can be approximated using the central difference quotient
	\begin{equation}
		\label{eq:approx_gradient}
		\nabla_{h,j} \mathcal{P}_h U_m(\bx)= \frac{\mathcal{P}_hU_m(\bx+he_j) - \mathcal{P}_hU_m(\bx-he_j)}{2h},  
	\end{equation}
	where $e_j$ is the unit vector in the $j$-th coordinate direction and the interested readers can refer to \cite{Lei2025}
	for more details.
	In terms of storage cost, only the grids $Q_k$ containing particles need to be stored, because $U_m(\bx)$ is approximated as $\mathcal{P}_h{U}_m(\bx) = 0$ in grids without particles and thus does not require storage. The dynamic grids generation based on particle locations (instead of full tensor-product grids) is why this method is termed the virtual uniform grid.

In short, SPM incorporates the nonlinear term into the relocation step, allowing particles to better capture the updated solution's features and maintain the desired distribution. This adjustment is expected to reduce variance and improve the accuracy, especially when dealing with strong nonlinearities, which has been validated in \cite{Lei2025} through numerical experiments. While the numerical convergence study about the parameters has been conducted previously in \cite{Lei2025}, we take the advection linear operator $\mathcal{L}=\bb\cdot\nabla$ as an example and establish the preliminary theoretical justification of SPM in the following theorem for the first time. 
\begin{theorem} 
	\label{thm:error_SPM}
	Under some regularity assumptions, when the linear operator has the form $\mathcal{L}=\bb\cdot\nabla$ and time stepsize $\tau\le A-1(A>1)$, the error between SPM solution $\tilde{U}_m$ and Lawson-Euler solution $U_m$ at time $t_m$ satisfies
	\begin{equation}
		\label{eq:error_SPM}
        \mathbb{E}\left[\|U_m - \tilde{U}_{m}\|_{2}^2\right] \le \me^{(B+1)T}\left(\mathcal{E}_{f}+\mathcal{E}_r(m)\right), 
    \end{equation}
	where $\|\cdot\|_2$ denotes the $L^2$-norm, $\mathcal{E}_{r}(m)$ represents the statistical error and spatial bias from the relocating up to time $t_m$, and $\mathcal{E}_{f}$ denotes the reconstruction error for the nonlinear term $f$, given by 
	\begin{equation}
		\begin{aligned}
			\mathcal{E}_{f} = AC_1h^2, \quad
			\mathcal{E}_{r}(m) = \sum_{j=0}^{m-1}
			\left(\frac{\mathbb{E}[Z_j^2]}{Nh^d}+C_0^2|U_j|_{H^1}^2 h^2\right),
		\end{aligned}
	\end{equation}
    with some constants $C_0=\frac{\sqrt{d}}{\pi}$, $C_1>0$ relying on the smoothness of $U_m$ and $f_m$, 
    and
    \begin{equation}\label{eq:B}
    B = \begin{cases}
    2AL_f^2(1+\frac{2}{h^2}), & \text{if $f$ depends on $\nabla u$}, \\
    2AL_f^2, & \text{if $f$ is independent of $\nabla u$}.
    \end{cases}
    \end{equation}
    Here $L_f$ is the Lipschitz constant of $f$ and $Z_j$ is the normalization constant defined as Eq.~\eqref{eq:relocating}. 
\end{theorem}

\begin{corollary}
    Under the same assumptions as Theorem \ref{thm:error_SPM}, the error between SPM solution $\tilde{U}_m$ and the exact solution $u$ satisfies
    \begin{equation}
        \begin{aligned}
            \mathbb{E}\left[\|u(\cdot, t_m)-\tilde{U}_{m}\|_{2}^2\right] \le& C_T\left(\tau^2+\mathcal{E}_f+\mathcal{E}_r(m)\right), 
        \end{aligned}
    \end{equation}
    with the constant $C_T\geq\me^{(B+1)T}$, where $B, \mathcal{E}_f$, $\mathcal{E}_r(m)$ are defined in Theorem \ref{thm:error_SPM}.
\end{corollary}

The error estimation indicates that SPM achieves a convergence rate of  $\mathcal{O}(\tau)$ in time, $\mathcal{O}(h)$ in space, and $\mathcal{O}(N^{-1/2})$ in the initial number of particles, as the existing numerical results have verified in \cite{Lei2025}. The regularity assumptions and detailed proof can be found in Section~\ref{sec:error_estimation}. 

\section{SPM-birth-death}\label{par bir}

SPM presents a comprehensive and effective framework with adaptive nature for high-dimensional problems \cite{Lei2025}. However, as Theorem~\ref{thm:error_SPM} reveals, the indiscriminate relocation mechanism at each time step $t_m$ accumulates both statistical error $\sqrt{Z_m^2/(Nh^d)}$ and spatial bias $C_0|U_m|_{H^1}h$, while incurring a computational cost of $\mathcal{O}(N)$ per time step. To address these issues, one straightforward idea is to perform resampling less frequently. Also, to avoid a gradual deviation of the particle distribution from the ideal distribution $p(\bx)\propto |U_m(\bx)|$ over time, an active birth mechanisms can be introduced to capture the small increment of order $\mathcal{O}(\tau)$ during the solution updation. Such idea motivates SPM-birth-death. 

SPM-birth-death adds an active particle birth-death mechanism on the basis of SPM \cite{Lei2025}. It inherits the main framework and design ideology of SPM, but differs in the specific particle evolution rules. High-dimensional nonlinear PDEs of the same form as Eq.~\eqref{cauchy problem} is considered. Also, consistent with SPM \cite{Lei2025}, SPM-birth-death uses the Lawson-Euler scheme, shown in Eq.~\eqref{eq:lawson}, to establish the evolution rules of $u(\bx,t)$ in the time direction. 

The fundamental difference between SPM and SPM-birth-death lies in the number of particles $N(t)$, which changes with time in SPM-birth-death, while remains constant in SPM. Same as SPM, each particle contains two parameters: the location $\bx_i \in \mathbb{R}^d$ and the weight $w_i\in \mathbb{R}$. All particles form a weighted point distribution
\begin{equation}\label{Xt}
	X_{t_m} := \frac{1}{N(0)} \sum_{i=1}^{N(t_m)} w_i(t_m) \delta_{\bx_i(t_m)}.
\end{equation}
where $w_i(t_m)$ and $\bx_i(t_m)$ denote the weight and location at time $t_m$, respectively. As a result, $X_{t_m}$ can approximate $U_m(\bx)$ in the weak sense
\begin{equation}\label{weak formulation}
		\langle \varphi, U_m\rangle
		\approx \langle \varphi,  X_{t_m} \rangle 
		= \frac{1}{N(0)} \sum_{i=1}^{N(t_m)} w_i(t_m) \varphi(\bx_i(t_m)),
\end{equation}
where the set of test functions is chosen as $\varphi(\bx) = \mone_{Q_{k}}(\bx), k = 1,2,\dots,K$, the same as SPM. Similarly, SPM-birth-death designs the evolution rules of the particle system $X_t$ based on the weak form of the Lawson-Euler scheme
\begin{equation}\label{eq:Lawson weak}
	\langle \varphi, U_{m+1} \rangle 
	= \langle \me^{\tau \mathcal{L}^*} \varphi, U_m + \tau  f(t_m,\bx,U_m,\vec{\nabla }U_m) \rangle. 
\end{equation}

\subsection{The birth mechanism}
We note that the particle system $X_{t_m}$ can already approximate $U_m(\bx)$ in the weak sense at time $t_m$. When the solution evolves from $U_m(\bx)$ to $U_{m+1}(\bx)$, it only adds an increment $\tau f\left(t_m, \boldsymbol{x}, U_m, \boldsymbol{\nabla} U_m\right)$ of order $\mathcal{O}(\tau)$, followed by the action of the operator $\me^{\tau \mathcal{L}^*}$.
Therefore, we can retain the original particle system that can approximate $U_m(\bx)$ at $t_m$, and then sample a small number of new particles of order $\mathcal{O}(\tau)$ from
\begin{equation}\label{C source}
	{p(\bx)}=\frac{|\tau  f(t_m,\bx,U_m,\vec{\nabla }U_m)|}{\int_{\mathbb{R}^d} |\tau  f(t_m,\bx,U_m,\vec{\nabla }U_m)| \D \bx},
\end{equation}
and add these new particles to the original particle system. 
Algorithm \ref{birdie alg} presents the detailed process of SPM-birth-death.

In Algorithm \ref{birdie alg}, $N(0)$ particles are sampled at the initial moment. At each subsequent time step, we sample $N_{\text{birth}}$ particles from Eq.~\eqref{C source} with $N_{\text{birth}}=N(0) \times \tau \int_{\mathbb{R}^d} | f(t_m,\bx,U_m,\vec{\nabla }U_m)| \D \bx$ (Line \ref{gen} in Algorithm \ref{birdie alg}), which reflects the magnitude of the nonlinear term and ensure the weak consistency of $L^1$ conservation. These new particles are sampled according to the distribution $\bx_i\sim {p(\bx)}=\frac{|\tau  f(t_m,\bx,U_m,\vec{\nabla }U_m)|}{\int_{\mathbb{R}^d} |\tau  f(t_m,\bx,U_m,\vec{\nabla }U_m)| \D \bx}$, and assigned the weight $w_i = \frac{|\tau  f(t_m,\bx_i,U_m,\vec{\nabla }U_m)|}{\tau  f(t_m,\bx_i,U_m,\vec{\nabla }U_m)}$ (Lines \ref{gen}-\ref{update w} in Algorithm \ref{birdie alg}). After sampling the new particles, we need to apply the operator $\me^{\tau \mathcal{L}^*}$. In addition, the updates of $\bx_i$ and $w_i$ require us to reconstruct the values of $U_m$ and $f(t_m,\bx,U_m,\vec{\nabla }U_m)$ from the particle point cloud. In the following, we present the specific implementation of these two aspects, which are consistent with those of SPM \cite{Lei2025}. 

\begin{itemize}
	\item The linear-operator-guided particle trajectory dynamics: Consider the linear equation,
	\begin{equation}\label{adjoint operator linear equation}
		\begin{aligned}
			\frac{\partial }{\partial t} g (\bx, t) &= \mathcal{L}^\ast g(\bx, t),\\  g(\bx, 0) &= \varphi(\bx),
		\end{aligned}
	\end{equation}
	and it has the formal solution $g(\bx,\tau) = \me^{\tau \mathcal{L}^*} \varphi(\bx)$, which reflects the action of $\mathcal{L}^*$. Here are some simple examples. If $\mathcal{L} = \vec{b}\cdot\vec{\nabla}$, $\vec{b}\in \mathbb{R}^d$ is a constant vector, then $\mathcal{L}^\ast = -\vec{b}\cdot\vec{\nabla}$ and the analytical solution of Eq.~\eqref{adjoint operator linear equation} is
	\[
	\me^{\tau \mathcal{L}^*} \varphi(\bx) = g(\vec{x},\tau) = \varphi(\vec{x}-\vec{b}\tau),
	\]
	which implies a kind of linear advection along the characteristic line,
	\begin{equation}\label{simu convec}
		\bx_i(t_{m+1}) = \bx_i(t_m) - \vec{b}\tau.
	\end{equation}
	If $\mathcal{L} = c\Delta$, $c$ is a constant, then $\mathcal{L}^\ast = c \Delta$, and the analytical solution of Eq.~\eqref{adjoint operator linear equation} has the closed form
	\[
	g(\vec{x},\tau) = \int_{\mathbb{R}^d} \frac{1}{(2\pi)^{d/2}\sqrt{|\vec{\Sigma}|}}  \me^{-\frac{1}{2} \vec{y}^T\vec{\vec{\Sigma}}^{-1}\vec{y}} \varphi(\vec{x}+\vec{y})\D \vec{y}, \quad \vec{\Sigma} = 2c \tau \vec{I},
	\]
	with $\vec{I}$ is the $d\times d$ identity matrix. It has a transparent probability interpretation $\mathrm{e}^{ c\tau \Delta } \varphi(\vec{x}) = \mathbb{E} \left[\varphi(\vec{x}+\vec{y})\right]$, where $\vec{y}$ obeying
	the normal distribution $\mathcal{N}(\vec{0}, \vec{\Sigma})$. To be more specific, it relates to the Brownian motion,
	\begin{equation}\label{simu diff}
		\vec{x}_i(t_{m+1}) = \vec{x}_i(t_m) + \vec{y}, \quad \vec{y} \sim \mathcal{N}(\vec{0}, \vec{\Sigma}).
	\end{equation}
	The above shows two examples of local linear operators. For the nonlocal linear operator $\mathcal{L}$, the particle motion rules can be designed through the Neumann series expansion \cite{rudin1991functional, kress2014linear}, and the interested readers
	can refer to  \cite{ShaoXiong2020, Lei2025} for more details.

	\item Constructing
	function values from the point cloud: The reconstruction of $U_m(\bx)$ inherits the VUG approach from SPM \cite{Lei2025}, as recalled in Section~\ref{sec:review_SPM}. That is, we derive the funtion and its derivative values using Eqs.~\eqref{eq weight sum approx u} and \eqref{eq:approx_gradient}.

\end{itemize}

Compared with SPM, birth mechanism reduces the computational cost of particle sampling and the errors induced by sampling. Note that at each time step, SPM samples $N(0)$ particles from
\begin{equation*}
	\frac{|U_{m}(\bx) + \tau  f(t_m,\bx,U_m,\vec{\nabla }U_m)|}{\int_{\mathbb{R}^d} |U_{m}(\bx) + \tau  f(t_m,\bx,U_m,\vec{\nabla }U_m)| \D \bx}.
\end{equation*}
This is unnecessary because the existing particle system can already approximate the part of $U_m(\bx)$. SPM-birth-death method only samples a small number of new particles from the $\mathcal{O}(\tau)$-order increment, which not only saves the computational cost but also reduces the additional errors introduced by sampling. 

\subsection{The death mechanism}
\par Now we have completed the introduction of the particle generation mechanism, which may cause the number of particles to grow rapidly with time during the simulation. Long-term simulations would be unfeasible without a corresponding particle annihilation mechanism to match it. Next, we discuss the particle annihilation mechanism in SPM-birth-death.

 
In the branching random walk algorithm for the Wigner equation \cite{ShaoXiong2020, XiongShao2020Overcoming}, the highly oscillatory nature of the Wigner function causes positive- and negative-weight particles to cluster densely. Consequently, particle annihilation in \cite{XiongShao2020Overcoming} operates through mutual cancellation of proximal positive/negative-weight particle pairs. For PDE solutions exhibiting weaker oscillations, particle systems typically concentrate positive-weight particles in solution-positive regions and negative-weight particles in solution-negative regions. Notably, PDEs with strictly non-negative solutions may contain no negative-weight particles. This inherent limitation restricts annihilation mechanisms relying on positive-negative cancellation. In \cite{Yan2015, Yan2016}, particle annihilation is implemented via resampling: Reconstruct the PDE solution from the $N(t_m)$-particle ensemble using the Fourier bases; release memory occupied by the original $N(t_m)$ particles;
sample $N(0)$ new particles from the Fourier-approximated solution; and annihilation occurs when $N(0) < N(t_m)$.

SPM-birth-death similarly reduces particle count through resampling with fewer particles. Unlike Fourier basis reconstruction—which employs tensor product structures exhibiting rapid dimensional scaling—our method utilizes the VUG approach \cite{Lei2025} for solution reconstruction. VUG capitalizes on particle adaptivity, achieving considerable storage reduction. The complete procedure remains detailed in Algorithm \ref{birdie alg} where $n_A \in (1, +\infty)$. This algorithm features the parameter $n_A$ to regulate particle population size. When the particle count $N(t_m)$ exceeds $N(0) \times n_A$, the particle annihilation mechanism is activated: The memory occupied by the original particle system is first released, and then $N(0)$ particles are sampled based on the piecewise constant approximation of the solution,
\begin{equation}\label{resample}
	\begin{aligned}
		\bx_i &\sim \frac{1}{\tilde{Z}_m}\sum_{k=1}^{K}\left|\tilde{U}_m(\bx)\right| \mone_{Q_{k}}(\bx), \quad i = 1,\dots,N(0),\\
		w_i &= \tilde{Z}_m \cdot \frac{\tilde{U}_m(\bx)}{\left|\tilde{U}_m(\bx)\right|}, \quad i = 1,\dots,N(0),
	\end{aligned}
\end{equation}
where $\tilde{Z}_m = \sum_{k=1}^{K}\left|\tilde{U}_m(\bx)\right| h^d$ is the normalization constant. 

In short, Algorithm \ref{birdie alg} provides a comprehensive description of SPM-birth-death. The moving and generation of particles are both guided by the weak form of the Lawson-Euler scheme, guaranteeing that the particle system consistently approximates the PDE solution in a weak sense throughout the simulation. Then to ensure the long-time evolution, the annihilation of particles is not directly derived from the weak form,  but is introduced as an auxiliary but efficient mechanism to control particle count under acceptable error and computational cost. The frequency of death mechanism can be flexibly  tuned by adjusting the parameter $n_A$. Here the physical interpretation is not as straightforward as particle methods for classical kinetic equations, because the particles do not represent physical particles any more, but serve as numerical carriers to approximate the solution.

Therefore, for the nonlinear PDE of the form~\eqref{cauchy problem}, SPM-birth-death is designed to possess the consistency and convergence as a numerical scheme, as Theorem~\ref{thm:error_SPM_birth_death} states. 
\begin{theorem}
	\label{thm:error_SPM_birth_death}
	Under some regularity assumptions, when the linear operator has the form $\mathcal{L}=\bb\cdot\nabla$ and time stepsize $\tau\le A-1(A>1)$, the error between SPM-birth-death solution $\tilde{U}_m$ and Lawson-Euler solution $U_m$ at time $t_m$ satisfies
	\begin{equation}
		\label{eq:error_SPM_birth_death}
		\mathbb{E}\left[\|U_m - \tilde{U}_{m}\|_{2}^2\right] \le \me^{(B+1)T}\left( \mathcal{E}_b + \mathcal{E}_f + \mathcal{E}_d \right)
	\end{equation}
	where $\mathcal{E}_b$ represents the statistical error introduced by the birth mechanism over time and $\mathcal{E}_d$ accounts for the additional statistical error and spatial bias from the death mechanism, given by
	\begin{equation}
		\begin{aligned}
			\mathcal{E}_b &= \frac{\mathbb{E}[\sum_{j=0}^{m-1}N_b(t_j)]}{N(0)^2h^d},\\
			\mathcal{E}_d &= (\mathbb{E}[R]+1)\max_{0\le k\le m}\left(\frac{\mathbb{E}[\tilde{Z}_k^2]}{N(0)h^d}+C_0^2|U_k|_{H^1}^2h^2\right).
		\end{aligned}
	\end{equation}
	Here $C_0, C_1, B$ and $\mathcal{E}_f$ are defined as Theorem~\ref{thm:error_SPM}, $N_b(t_m)$ denotes the number of particles born at the instant $t_m$, $\tilde{Z}_m$ is the normalization constant of $\tilde{U}_m$, and $R$ is the number of annihilation steps, satisfying
    \begin{equation}\label{eq:R}
        \mathbb{E}[R] \le \min\left\{\frac{\mathbb{E}[\sum_{j=0}^{\frac{T}{\tau}-1}N_b(t_j)]}{(n_A-1)N(0)}, \frac{T}{\tau}\right\}. 
    \end{equation}
\end{theorem}

Combining Theorem \ref{thm:error_SPM_birth_death} with the first-order accuracy of Lawson-Euler scheme together yields immediately the overall error estimation between the SPM-birth-death solution and the exact solution.

\begin{corollary}
    Under the same assumptions as Theorem \ref{thm:error_SPM_birth_death}, the error between SPM-birth-death solution $\tilde{U}_m$ and the exact solution $u$ satisfies
    \begin{equation}
        \begin{aligned}
            \mathbb{E}\left[\|u(\cdot, t_m)-\tilde{U}_{m}\|_{2}^2\right] \le& C_T\left(\tau^2+\mathcal{E}_b+\mathcal{E}_f+\mathcal{E}_d\right), 
        \end{aligned}
    \end{equation}
    with the constant $C_T\geq\me^{(B+1)T}$, where $B, \mathcal{E}_b$, $\mathcal{E}_f$ and $\mathcal{E}_d$ are defined in Theorem \ref{thm:error_SPM_birth_death}.
\end{corollary}

The theoretical result above indicates that SPM-birth-death achieves the same convergence order of $\mathcal{O}(\tau)$ in time, $\mathcal{O}(h)$ in space and $\mathcal{O}(N(0)^{-1/2})$ in the initial particle number as SPM. The statistical error has the order $\mathcal{O}((N(0)h^d)^{-1/2})$ contributed by the initial $N(0)$ particles, which is similar to the classical Monte Carlo method and the same trick of (conditional) independence is also employed in the later proof. The $h^{d}$ in the denominator comes from the classical histogram density estimation error analysis, which indicates that smaller cell size $h$ leads to overfitting and larger variance. Consequently, $h$ would not have been chosen too small in practice, so that the exponential factor $B$ including $1+\frac{2}{h^2}$ when $f$ depends on $\nabla u$ would not be too large, see Eq.~\eqref{eq:B}. In fact, the $1+\frac{2}{h^2}$ factor in $B$ is essential because it stems from the gradient estimation error when using the $L^2$-norm. 


Our theoretical analysis also reveals that SPM-birth-death, with suitable $n_A$ and thus less frequent resampling, can reduce the total error compared with SPM. Note that the differences of error estimations lie in the resampling error $\mathcal{E}_r$ introduced by the relocating step in SPM and $\mathcal{E}_d$ introduced by the death mechanism in SPM-birth-death, as well as statistical error $\mathcal{E}_b$ caused by the birth mechanism. On one hand, the magnitude of both resampling errors is related to the number of resampling steps. That is, the relocating error $\mathcal{E}_r$ in SPM accumulates at every time step, while the death error $\mathcal{E}_d$ in SPM-birth-death only accumulates when the particle number exceeds the threshold. On the other hand, the additional $\mathcal{E}_b$ introduced by SPM-birth-death is of order $\mathcal{O}(\frac{1}{N(0)h^d})$ and is negligible compared with $\mathcal{E}_r$ when $\tau$ is small. This comparison elucidates the indiscriminate annihilation in SPM may lead to unnecessary error accumulation, while SPM-birth-death can effectively mitigate this by adjusting the frequency of resampling.

The regularity assumptions and detailed proof of Theorem~\ref{thm:error_SPM_birth_death} are deferred to Section~\ref{sec:error_estimation}, followed by numerical experiments in Section~\ref{par bir numer} to validate the theoretical results as well as to demonstrate the efficiency of SPM-birth-death.

\begin{algorithm}[htb]
	\caption{ Stochastic particle method with birth-death dynamics. } 
	\label{birdie alg}
	\hspace*{0.02in} {\bf Input:}
	The side length of hypercube $h$, the final time $T$, the time step $\tau$, the particle growth threshold $n_A$, the initial data $u_0(\bx)$, 
	the linear operator $\mathcal{L}$, the nonlinear term {$f(t,\bx,u,\vec{\nabla }u)$} and the initial partical number $N(0)$.\\
	\hspace*{0.02in} {\bf Output:} The numerical solution to nonlinear PDE \eqref{cauchy problem} at the final time.
	\begin{algorithmic}[1]
		\For{$m=0:\frac{T-\tau}{\tau}$} 
		\If{$N(t_m) > N(0) \times n_A$} \label{ann begin}
		\State Release the original $N(t_m)$ particles.
		\State Sample $N(0)$ particles based on the piecewise constant approximation of $U_m$ according to \eqref{resample};
		\State $N(t_m)\leftarrow N(0)$;
		\EndIf \label{ann end}
		\If{$m = 0$}
		\State Sample $\bx_i$ from $ p(\bx)=\frac{|U_0(\bx)|}{\int_{\mathbb{R}^d} |U_0(\bx)| \D \bx},\quad i = 1,\dots,N(0)$;
		\State $w_i\leftarrow\frac{U_0(\bx_i)}{{p(\bx_i)}},\quad i = 1,\dots,N(0)$;
		\EndIf
		\State Compute the particle number to be generated $N_{\text{birth}} \leftarrow N(0) \times \tau \int_{\mathbb{R}^d} | f(t_m,\bx,U_m,\vec{\nabla }U_m)| \D \bx$; \label{gen}
		\State Sample new particles $\bx_i$ from ${p(\bx)}=\frac{|\tau  f(t_m,\bx,U_m,\vec{\nabla }U_m)|}{\int_{\mathbb{R}^d} |\tau  f(t_m,\bx,U_m,\vec{\nabla }U_m)| \D \bx}, \quad i = N(t_m)+1,\dots,N(t_m)+N_{\text{birth}}$; \label{update x}
		\State $w_i\leftarrow \frac{|\tau  f(t_m,\bx_i,U_m,\vec{\nabla }U_m)|}{\tau  f(t_m,\bx_i,U_m,\vec{\nabla }U_m)}, \quad i = N(t_m)+1,\dots,N(t_m)+N_{\text{birth}}$; \label{update w}
		\State $N(t_{m+1})\leftarrow N(t_m)+N_{\text{birth}}$;
		\State Move particle according to $\mathcal{L}^*$;
		\State Obtain the piecewise constant approximation of $U_{m+1}$ and $f(t_{m+1},\bx,U_{m+1},\vec{\nabla }U_{m+1})$ based on the virtual uniform grid;
		\EndFor
	\end{algorithmic}
\end{algorithm}

\subsection{Complexity and parameter optimization}
Here we give a brief discussion on the computational complexity of SPM-birth-death, which guides the practical selection of parameter $n_A$. We will also compare the complexity with that of SPM in \cite{Lei2025}, which helps to illustrate the efficiency of the introduced active birth-death mechanism. For time step $t_m$, the computational cost mainly comes from four parts: 
\begin{enumerate}[(a)]
  \item the resampling of particles has a cost of $\mathcal{O}(N(0))$ if the death mechanism occurs; 
  \item the birth process has a cost of $\mathcal{O}(\mathbb{E}[N_b(t_m)])=\mathcal{O}(\tau N(0))$; 
  \item the evolution of particles according to the linear part has a cost of $\mathcal{O}(N(t_m))$; 
  \item the piecewise constant reconstruction  has a cost of $\mathcal{O}(N(t_m)\log N(t_m))$ \cite{Lei2025}. 
\end{enumerate}
Therefore, the total computational cost can be estimated as
\begin{equation}
    \mathcal{O}\left(RN(0)+ \sum_{m=0}^{M-1} N(t_m) \log N(t_m) \right). 
\end{equation}
Note that between two consecutive death mechanisms at time steps $m_r$ and $m_{r+1}$, the sum of particle number can be bounded as 
\begin{equation*}
    \sum_{m=m_r}^{m_{r+1}-1} N(t_m) \log N(t_m) \le (m_{r+1}-m_r) n_A N(0) \log(n_A N(0)) + \mathcal{O}(\tau),
\end{equation*}
where the last term comes from the birth process at $t_{m_{r+1}-1}$. If we are allowed to neglect the logarithm term, 
according to Eq.~\eqref{eq:R}, 
the expectation of total computational cost can be further bounded as
\begin{equation}
    \mathcal{O}\left(\frac{T}{\tau} n_A N(0) + \mathbb{E}[R] N(0)\right) = \mathcal{O}\left(\frac{T}{\tau} n_A N(0) + \min\left\{\frac{\mathbb{E}[\sum_{j=0}^{M-1}N_b(t_j)]}{(n_A-1)N(0)}, \frac{T}{\tau}\right\} N(0)\right), 
\end{equation}
Therefore, to minimize the computational cost, we can choose the prescribed threshold for death  
\begin{equation}
    n_A^* = \mathcal{O}(\sqrt{\frac{\mathbb{E}[\tau\sum_{j=0}^{M-1}N_b(t_j)]}{T N(0)}}) +1 = \mathcal{O}(\sqrt{\frac{\tau}{T}}) + 1, 
\end{equation}
which implies that when $R=\mathcal{O}(\sqrt{T/\tau})$, the total computational cost can be minimized. The optimal selection of $n^*_A$ balances the cost of the death mechanism with that of the walk and reconstruction process, ensuring that neither of them dominates the total computational cost. In addition, Theorem~\ref{thm:error_SPM_birth_death} demonstrates that the error bound scales as $\mathcal{O}(\sqrt{1/(n_A-1)})$ with respect to $n_A$. Therefore, the main influence of $n_A$ is to control the frequency of annihilations, and thus balance the trade-off between the accuracy and the computational cost. We would like to point out that the optimal choice of $n_A$ may vary  and depend on the specific problem as well as the desired accuracy level.

Note that in SPM, the relocating mechanism is performed at every time step, which can be viewed as the special case of SPM-birth-death with a small $n_A$: $(n_A-1)N(0)<N_{birth}$, since the cost of the birth process is negligible compared with the cost of the relocating mechanism: $\mathcal{O}(N(0))$. Hence, the computational cost of SPM (corresponding to the small $n_A$) is larger than that of SPM-birth-death method with $n_A$ close to the optimal choice of $n^*_A$. 
Even worse, the indiscriminate annihilations in SPM also lead to larger error accumulation over time,
since every annihilation step introduces additional bias and variance, but 
SPM-birth-death only performs resampling when needed,
which can be also clearly seen by comparing Theorem~\ref{thm:error_SPM} with Theorem~\ref{thm:error_SPM_birth_death}.
Consequently, to achieve the same accuracy level, SPM-birth-death may require much lower computational complexity than SPM, particularly in high-dimensional regimes or long-time integrations.

\section{Numerical analysis}
\label{sec:error_estimation}

\begin{figure}[hbtp]
    \centering
    \includegraphics[width=.6\textwidth]{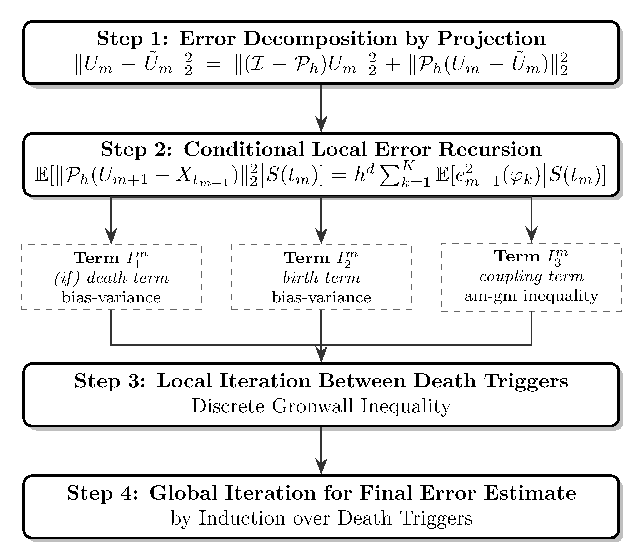}
    \caption{Flowchart for proving Theorem~\ref{thm:error_SPM_birth_death}, an error estimator of SPM-birth-death.}
    \label{fig:flowchartproof}
\end{figure}

This section presents the proof of Theorems~\ref{thm:error_SPM} and \ref{thm:error_SPM_birth_death}, estimating the error of SPM and SPM-birth-death, respectively. Section~\ref{subsec:assumptions_lemmas} introduces the necessary assumptions, settings and lemmas, followed by the proof in Section~\ref{subsec:proof_SPM_birth_death}. The overall proof structure of Theorem~\ref{thm:error_SPM_birth_death} is summarized in Figure~\ref{fig:flowchartproof}, while Theorem~\ref{thm:error_SPM} can be viewed as a special case of Theorem~\ref{thm:error_SPM_birth_death}, obtained by extending the proof framework of SPM-birth-death. 

\subsection{Assumptions, settings and lemmas}
\label{subsec:assumptions_lemmas}

To proceed with the theoretical analysis, we introduce the following assumptions. The first two assumptions about regularity are standard in the context of nonlinear parabolic equation, while the last one is technical and imposed primarily to simplify the proofs, which we believe does not compromise the generality of the method significantly.

\begin{assumption}
    \label{ass:nonlinear}
    The nonlinear term $f(t,\bx, u, \nabla u)$ satisfies the Lipschitz condition with respect to $u$ and $\nabla u$, i.e., there exists a constant $L_f>0$ such that for any $t\in[0,T]$, $\bx\in\mathbb{R}^d$, $u_1,u_2\in\mathbb{R}$ and $v_1, v_2\in\mathbb{R}^d$, we have
    \begin{equation}
        |f(t,\bx,u_1,v_1) - f(t,\bx,u_2,v_2)| \le L_f \sqrt{|u_1 - u_2|^2 + \|v_1 - v_2\|^2}.
    \end{equation}
\end{assumption}
\begin{assumption}
    \label{ass:sol}
    The exact solution $u$ to Eq.~\eqref{cauchy problem} possesses sufficient regularity and decays rapidly at infinity. Consequently, we assume that the Lawson-Euler solution $U_m$ and the nonlinear term $f_m$ inherit these properties, maintaining sufficient smoothness and fast decay rates.
\end{assumption}
\begin{assumption}
    \label{ass:linear}
    The linear operator has the form of $\mathcal{L} = \bb\cdot\nabla$, where $\bb\in\mathbb{R}^d$ is a constant vector, then the particle trajectory is deterministic. 
\end{assumption}


Before presenting the proof,  we first make some preparations. Let 
\begin{align*}
S(t_m) :=& \{\bx_i(t_m),w_i(t_m)\}_{i=1}^{N(t_m)}, \\
\tilde{f}_m :=& f(t_m,\bx,\tilde{U}_m,\nabla_h \tilde{U}_m), \\
f_m: =& f(t_m,\bx,U_m, \nabla U_m).
 \end{align*}
Conditional on the particle system $S(t_m)$, we have known $X_{t_m}$ and $\tilde{U}_m$. After one step of SPM-birth-death, we obtain the new particle system $S(t_{m+1})$ and the coorresponding weak approximation
\begin{equation}
    \label{eq:algorithm1step}
    X_{t_{m+1}} = \me^{\tau\mathcal{L}}(X_{t_m,o}+X_{t_m,b}), 
\end{equation}
where $X_{t_m,o}$ and $X_{t_m,b}$ denote the original or resampled particles after the death mechanism and the new born particles after the birth mechanism, respectively. More details are as follows.
\begin{itemize}
    \item Death mechanism: $X_{t_m,o}$ can be represented as
    \begin{equation}
        X_{t_m,o} = \frac{1}{N(0)}\sum_{i=1}^{N_o(t_m)} w_{i,o}(t_m) \delta_{\bx_{i,o}(t_m)}.
    \end{equation}
    If the death mechanism does not occur at time $t_m$, then $X_{t_m,o}=X_{t_m}$. Otherwise, the resampling procedure in Eq.~\eqref{resample} gives the distribution of particle locations and weights. That is, the location of each particle $\bx_{i,o}(t_m)$ in $Q_k$ is uniformly distributed and has the following variance in the weak sense
    \begin{equation}
      \label{eq:var_location}
      \Var_{\bx\sim \varphi_k}(\varphi(\bx)) = \frac{1}{|Q_k|}\int_{Q_k}(\varphi(\bx))^2\D \bx - \left(\frac{1}{|Q_k|} \int_{Q_k} \varphi(\bx)\right)^2. 
    \end{equation}
    Meanwhile, the number of particles in a cell is a random variable following the binomial distribution 
    \begin{equation}\label{eq:binomial}
        |\{i\big| \bx_{i,o}(t_m)\in Q_k\}|=\begin{cases}
           [N_{o,k}(t_m)]+1, & \text{with probability } \{N_{o,k}(t_m)\},\\
            [N_{o,k}(t_m)], & \text{with probability } 1-\{N_{o,k}(t_m)\},
        \end{cases}
    \end{equation}
    where $N_{o,k}(t_m)=h^d\frac{|\tilde{U}_m|(Q_k)}{\tilde{Z}_m}N(0)$ is the conditional expectation of resampling particle numbers in the cell $Q_k$, $[\cdot]$ and $\{\cdot\}$ denote the integer part and fractional part, respectively. Thus the total weight in each cell $Q_k$ has the conditional expectation and variance
    \begin{equation}
      \label{eq:exp_var_weights}
        \mathbb{E}\left[\sum_{\bx_{i,o}(t_m)\in Q_k} w_{i,o} \Big| S(t_m)\right] = \sum_{\bx_i(t_m)\in Q_k} w_i,\quad        \Var\left[\sum_{\bx_{i,o}(t_m)\in Q_k} w_{i,o} \Big| S(t_m)\right] \le \frac{1}{4} \tilde{Z}_m^2.
    \end{equation}

    \item Birth mechanism: $X_{t_m,b}$ can be represented as
    \begin{equation}
        X_{t_m,b} = \frac{1}{N(0)}\sum_{i=1}^{N_b(t_m)} w_{i,b}(t_m) \delta_{\bx_{i,b}(t_m)},
    \end{equation}
    where the locations $\{\bx_{i,b}(t_m)\}_{i=1}^{N_b(t_m)}$ are i.i.d. samples from the distribution, i.e., 
    \begin{equation}
        \bx_{i,b}(t_m)\sim \frac{|\tau \mathcal{P}_h\tilde{f}_m(\bx)|}{\int_{\mathbb{R}^d} |\tau \mathcal{P}_h\tilde{f}_m(\bx)| \D \bx},\quad w_{i,b}(t_m) = \frac{\tau \mathcal{P}_h\tilde{f}_m(\bx_{i,b}(t_m))}{|\tau \mathcal{P}_h\tilde{f}_m(\bx_{i,b}(t_m))|}. 
    \end{equation}
    The distribution above is designed to guarantee the unbiased conditional expectation of $X_{t_m,b}$ in the weak sense, that is, for any test function $\varphi$,
    \begin{equation}
      \label{eq:birth_exp}
      \mathbb{E}[\langle\varphi,X_{t_m,b}\rangle|S(t_m)] = \langle \varphi, \tau \mathcal{P}_h\tilde{f}_m \rangle.
    \end{equation}
    And the number of newly born particle $N_b(t_m)$ is a random variable following the similar binomial distribution~\eqref{eq:binomial}, and satisfying
    \begin{equation}
      \mathbb{E}[N_b(t_m)|S(t_m)]=N_{birth}(t_m)=N(0)\times\int_{\mathbb{R}^d} |\tau \mathcal{P}_h\tilde{f}_m(\bx)| \D \bx, 
    \end{equation}
    to ensure the total weight of the born particles matching the increment $\tau \tilde{f}_m$. 
    
    \item Walk mechanism: for the particle trajectory part, denote $\varphi^t=\me^{\mathcal{L}^*t}\varphi$ as the evolution of the test function $\varphi$. The analytical expression of $\varphi^t$ can be derived by characteristic lines under Assumption \ref{ass:linear}, i.e., 
\begin{equation}
    \label{eq:analytical_sol}
    \varphi^t(\bx) = \varphi(\bx-\bb t). 
\end{equation}
\end{itemize}

Considering the procedure composed of three mechanisms above, we can readily figure out that 
the particle system has the Markovian property, that is, $S(t_{m+1})$ depends only on $S(t_m)$,
although such property may not hold for each single particle. 

We further need the following lemmas about the norm estimation of the evolved test functions. 
\begin{lemma}
    \label{lem:Linfupper}
    Under Assumption \ref{ass:linear}, given any test function $\varphi_k\in L^\infty(\mathbb{R}^d)(k=1,2,\dots)$ with pairwise disjoint support sets, and time $t\geq 0$, the following inequality holds, 
    \begin{equation}
        \|\sum_k(\varphi_k^t)^2\|_{\infty}\le \sup_k\|\varphi_k\|_\infty^2.
    \end{equation}
\end{lemma}
\begin{proof}[Proof of Lemma \ref{lem:Linfupper}]
    Note that the support sets of $\{\varphi_k\}$ are pairwise disjoint, then we get
    \begin{equation*}
        \|\sum_k(\varphi_k^t)^2\|_\infty = \|\sum_k \me^{\mathcal{L}^*t}(\varphi_k)^2\|_\infty 
        = \|\me^{\mathcal{L}^*t}(\sum_k(\varphi_k)^2)\|_\infty = \|\sum_k (\varphi_k)^2\|_\infty\le \sup_k \|\varphi_k\|_\infty^2.
    \end{equation*}
\end{proof}

\begin{lemma}
    \label{lem:L2preserve}
    Under Assumption \ref{ass:linear}, given any test function $\varphi\in L^2(\mathbb{R}^d)$, the $L^2$-norm is preserved over time, i.e,
    \begin{equation}
        \|\varphi^t\|_2 = \|\varphi\|_2.
    \end{equation}
\end{lemma}

To estimate piecewise constant projection error, we need the Poincar\'{e} inequality \cite{Payne1960}. 
\begin{lemma}
    \label{lem:poincare}
    For any function $v \in H^1(\mathbb{R}^d)$, the following Poincar\'{e} inequality holds:
    \begin{equation}\label{eq:P}
        \|v - \mathcal{P}_h v\|^2_2 \le C_0^2 |v|_{H^1}^2 h^2,
    \end{equation}
    where $\mathcal{P}_h$ is the piecewise constant projection operator in \eqref{eq weight sum approx u} and the optimal coefficient is $C_0=\frac{\sqrt{d}}{\pi}$.
\end{lemma}

To compute the cumulative error, we need the following discrete Gronwall inequality.
\begin{lemma}
    \label{lem:gronwall}
    For any $a_k\geq 0$, $B\geq 0$, $c_k\geq 0$, $\tau >0$, if the inequality $a_k \leq c_k + B\tau\sum_{m=0}^{k-1} a_m$ holds, then we have
    \begin{equation}
        a_k\le \me^{B k \tau}\max_{0\le m\le k}c_m. 
    \end{equation}
\end{lemma}
\begin{proof}[Proof of Lemma \ref{lem:gronwall}]
    Let $S_k=\sum_{m=0}^k a_m$, we have
    \begin{equation*}
        S_k \le c_k + S_{k-1}(1+\tau B), 
    \end{equation*}
    which can be iterated and rewritten as
    \begin{equation*}
        \begin{aligned}
            S_k &\le& \sum_{m=0}^{k} c_m (1+\tau B)^{k-m}\le \sum_{m=0}^k (1+\tau B)^{k-m} \max_{0\le m \le k} c_m =\frac{(1+\tau B)^{k+1}-1}{\tau B} \max_{0\le m \le k} c_m.
        \end{aligned}
    \end{equation*}
    Thus we obtain
    \begin{equation*}
        a_k \le c_k + B\tau S_{k-1} \le (1+\tau B)^k \max_{0\le m \le k} c_m \le \me^{Bk\tau} \max_{0\le m \le k} c_m.
    \end{equation*}
\end{proof}

\subsection{Proof of theorems}\label{subsec:proof_SPM_birth_death}

Following the flowchart in Figure~\ref{fig:flowchartproof}, we begin the proof of Theorem \ref{thm:error_SPM_birth_death}.

\begin{proof}[Proof of Theorem \ref{thm:error_SPM_birth_death}]
    Denote the error of the weak formulation~\eqref{weak formulation} by
    \begin{equation}
        e_m^2(\varphi) := (\langle\varphi, U_m-X_{t_m}\rangle)^2, 
    \end{equation}
    where $\varphi$ is a test function. Let $\varphi_k=\frac{1}{|Q_k|}\mone_{Q_k}$, and the projection of function $v$ can be represented as the linear combination of the basis $\{\varphi_k\}$, i.e.,
    \begin{equation}
        \mathcal{P}_h v = \sum_{k=1}^K h^d \langle \varphi_{k}, v\rangle \varphi_{k}.
    \end{equation}
    Denote the identity $\mathcal{I}v=v$, then we have, $\forall v\in H^1(\mathbb{R}^d)$, 
    \begin{align}
       h^d\sum_{k=1}^K (\langle \varphi_k,v \rangle)^2 = \|\mathcal{P}_h v\|_2^2&=\|v\|_2^2-\|(\mathcal{I}-\mathcal{P}_h)v\|_2^2 \label{eq:projection_identity} \\
        &\geq\|v\|_2^2-C_0^2|v|_{H^1}^2 h^2, \label{eq:projection_ieq}
    \end{align}
    where $C_0=\frac{\sqrt{d}}{\pi}$ is the absolute constant from Lemma \ref{lem:poincare}. Hence, we obtain the error decomposition by projection as Step 1 in Figure~\ref{fig:flowchartproof},
    \begin{equation}
        \label{eq:L2_error}
        \|\mathcal{P}_h({U}_m-\tilde{U}_m)\|_2^2 = h^d\sum_{k=1}^K e_m^2(\varphi_k), \quad
        \|{U}_m - \tilde{U}_m\|_2^2 \le h^d\sum_{k=1}^K e_m^2(\varphi_k) + C_0^2|{U}_m|_{H^1}^2 h^2. 
    \end{equation}
    
    Then we turn to Step 2 in the flowchart, the recursion formula for the conditional expectation of $e_m^2(\varphi)$. Following one step  procedure of SPM-birth-death as Eq.~\eqref{eq:algorithm1step} and the Lawson-Euler scheme in Eq.~\eqref{eq:lawson}, we have
    \begin{align}
                \mathbb{E}[e_{m+1}^2(\varphi)\big|S(t_m)] =& \mathbb{E}\left[\left(\langle \varphi^\tau, U_{m}-X_{t_m,o}\rangle+\langle\varphi^\tau, \tau f_m-X_{t_m,b}\rangle \right)^2\big|S(t_m)\right]\nonumber\\
            =:&I^m_1(\varphi^\tau)+I^m_2(\varphi^\tau)+I^m_3(\varphi^\tau). \label{eq:I123}
            \end{align}

\begin{itemize}
    \item {\it As for $I_1^m(\varphi^\tau)=\mathbb{E}[(\langle \varphi^\tau, {U}_{m}-X_{t_m,o}\rangle)^2\big|S(t_m)]$}. The term is related to the death mechanism. If the death mechanism does not occur at time $t_m$, then $X_{t_m,o}=X_{t_m}$ and we have
    \begin{equation}
        \label{eq:I1_error_no_death}
        I^m_1(\varphi^\tau) = \mathbb{E}\left[\left(\langle \varphi^\tau, {U}_{m}-X_{t_m}\rangle\right)^2\big|S(t_m)\right] = \left(\langle \varphi^\tau, {U}_{m}-X_{t_m}\rangle\right)^2=e_m^2(\varphi^\tau).
    \end{equation}
    Otherwise, it can be decomposed into spatial bias and resampling variance as
    \begin{equation}\label{eq:b-v}
        I^m_1(\varphi^\tau) =\left(\langle \varphi^\tau, {U}_{m}-\tilde{U}_m\rangle\right)^2 + \mathbb{E}\left[\left(\langle \varphi^\tau, \tilde{U}_m - X_{t_m,o}\rangle\right)^2\big|S(t_m)\right]:=I^m_{1,1}(\varphi^\tau)+I^m_{1,2}(\varphi^\tau).
    \end{equation}
    For the spatial bias part, taking $\varphi=\varphi_k$ for $I_{1,1}^m(\varphi^\tau)$, summing over $k=1,2,\dots,K$,
    and using the projection identity~\eqref{eq:projection_identity} and Lemma \ref{lem:L2preserve} lead to
    \begin{equation}
        \begin{aligned}
            \label{eq:I11_error}
            h^d\sum_{k=1}^K I^m_{1,1}(\varphi_k^\tau) &=  h^d \sum_{k=1}^K \left(\langle \varphi_k, \me^{\tau \mathcal{L}}({U}_{m}-\tilde{U}_m)\rangle\right)^2
            = \|\mathcal{P}_h \me^{\tau \mathcal{L}}({U}_{m}-\tilde{U}_m)\|_2^2\le \|{U}_{m}-\tilde{U}_m\|_2^2.
        \end{aligned}
    \end{equation}
    For the reconstruction variance part, we seperate the randomness of locations and weights by the expresion of piecewise constant approximation~\eqref{eq weight sum approx u} and obtain
    \begin{equation*}
        \begin{aligned}
            \langle \varphi^\tau, X_{t_m, o} - \tilde{U}_m\rangle =&\sum_{k=1}^K \sum_{\bx_{i,o}(t_m) \in Q_k}  \frac{w_{i,o}}{N(0)}\left(\varphi^\tau(\bx_{i,o}(t_m)) - \frac{1}{|Q_k|} \int_{Q_k} \varphi^\tau(\vb*{x}) d\vb*{x}\right)\\
            & + \sum_{k=1}^K \left(\frac{1}{|Q_k|} \int_{Q_k} \varphi^\tau(\vb*{x}) d\vb*{x}\right) \frac{1}{N(0)}\left(\sum_{\bx_{i,o}(t_m)\in Q_k} w_{i,o} - \sum_{\bx_i(t_m)\in Q_k} w_i\right).
        \end{aligned}
    \end{equation*}
    In consequence, using the independence of particle locations, weights and number, we have 
    \begin{align}
                I^m_{1,2}(\varphi^\tau)=& \mathbb{E}\left[\left(\langle \varphi^\tau, \tilde{U}_m - X_{t_m,o}\rangle\right)^2\big|S(t_m)\right] \nonumber \\
            =& \frac{1}{N(0)^2}\sum_{k=1}^K\mathbb{E}\left[\left(\sum_{\bx_{i,o}(t_m) \in Q_k} (w_{i,o})^2\right)\big|S(t_m)\right]\Var_{\bx\sim \varphi_k}\left(\varphi^\tau(x)\right) \nonumber \\
             &+ \frac{1}{N(0)^2}\sum_{k=1}^K \left(\frac{1}{|Q_k|} \int_{Q_k} \varphi^\tau(\vb*{x}) d\vb*{x}\right)^2\Var\left[\left(\sum_{\bx_{i,o}(t_m)\in Q_k} w_{i,o}\right)\big|S(t_m)\right] \nonumber \\
            \le& \frac{1}{N(0)^2}\sum_{k=1}^K\mathbb{E}\left[\left(\sum_{\bx_{i,o}(t_m) \in Q_k} (w_{i,o})^2\right)\big|S(t_m)\right]\frac{1}{|Q_k|}\int_{Q_k}(\varphi^\tau(\vb*{x}))^2\dd \vb*{x}\nonumber \\
            =& \frac{\tilde{Z}_m}{N(0)}\langle (\varphi^\tau)^2, |\tilde{U}_m|\rangle,  \label{eq:I12}
    \end{align}
    where the inequality holds due to Eqs.~\eqref{eq:var_location} and \eqref{eq:exp_var_weights}. Therefore, take $\varphi=\varphi_k$, sum over $k=1,2,\dots,K$, and employ Lemma \ref{lem:Linfupper} to obtain
    \begin{equation}
        \begin{aligned}
            \label{eq:I12_error}
            h^d\sum_{k=1}^K I^m_{1,2}(\varphi_k^\tau) \le \frac{\tilde{Z}_m}{N(0)}\langle h^d\sum_{k=1}^K (\varphi_k^\tau)^2, |\tilde{U}_m|\rangle
            \le \frac{\tilde{Z}_m}{N(0)}h^d\|\sum_{k=1}^K (\varphi_k^\tau)^2\|_{\infty}\|\tilde{U}_m\|_1\le \frac{\tilde{Z}_m^2}{N(0) h^d}.
        \end{aligned}
    \end{equation}
    In particular, $I_1^0$ is just the initial statistical error and we can similarly bound it as Eq.~\eqref{eq:I12_error}.

    \item {\it As for $I_2^m(\varphi^\tau)=\mathbb{E}[(\langle\varphi^\tau, \tau {f}_m-X_{t_m,b}\rangle )^2\big|S(t_m)]$}. The term is related to the birth mechanism. By a similar bias-variance decomposition as Eq.~\eqref{eq:b-v} for $I_1^m$, we have
    \begin{equation}\label{eq:I2}
        \begin{aligned}
            I^m_2 = \left(\langle \varphi^\tau, \tau {f}_m - \tau \mathcal{P}_h\tilde{f}_m\rangle\right)^2 + \mathbb{E}\left[\left(\langle \varphi^\tau, \tau \mathcal{P}_h\tilde{f}_m - X_{t_m,b}\rangle\right)^2\big|S(t_m)\right]=:&I^m_{2,1}+I^m_{2,2}. 
        \end{aligned}
    \end{equation}
       
       Under the Lipschitz continuity of nonlinear term $f$ in Assumption \ref{ass:nonlinear}, we can further bound the error terms appropriately, i.e., 
    \begin{equation}\label{eq:ff}
        \begin{aligned}
            \|{f}_m-\tilde{f}_m\|_2^2\le L_f^2\left(\|{U}_m-\tilde{U}_m\|_2^2+\|\nabla{U}_m-\nabla_h\tilde{U}_m\|_2^2\right).
        \end{aligned}
    \end{equation}
    Note that the central difference approximation $\nabla_h$ to the gradient satisfies
    \begin{equation*}
        \|\nabla U_m-\mathcal{P}_h \nabla_hU_m\|_2 \le \|(\mathcal{I}-\mathcal{P}_h)\nabla U_m\|_2 + \|\mathcal{P}_h(\nabla U_m - \nabla_hU_m)\|_2\le C_u h,
    \end{equation*}
    for some constant $C_u>0$ depending on the regularity of $U_m$, where the last inequality holds because of the Poincar\'e inequality~\eqref{eq:P} and the consistency of the central difference scheme. Since $\nabla_h$ and $\mathcal{P}_h$ are commutative, we have
    \begin{equation}\label{eq:gUU}
        \|\nabla U_m-\nabla_h\tilde{U}_m\|_2\le \|\nabla U_m-\mathcal{P}_h \nabla_hU_m\|_2+\|\nabla_h(\mathcal{P}_h{U}_m-\tilde{U}_m)\|_2\le C_u h + \frac{1}{h}\|\mathcal{P}_h{U}_m-\tilde{U}_m\|_2.
    \end{equation}
    Also, we have 
    \begin{equation}\label{eq:UU}
    \|U_m-\tilde{U}_m\|_2^2\le \|\mathcal{P}_hU_m-\tilde{U}_m\|_2^2+C_0^2 |U_m|_{H^1}^2h^2. 
    \end{equation}
    
    Hence, for the spatial bias part of Eq.~\eqref{eq:I2}, 
    we have
    \begin{align}
        h^d\sum_{k=1}^K I^m_{2,1}(\varphi_k^\tau) \le& \tau^2 \|{f}_m-\mathcal{P}_h\tilde{f}_m\|_2^2 \nonumber \\
        \le& 2\tau^2 (\|(\mathcal{I}-\mathcal{P}_h){f}_m\|_2^2+\|\mathcal{P}_h({f}_m-\tilde{f}_m)\|_2^2)\nonumber\\
        \le& 2\tau^2(C_0^2|f_m|_{H^1}^2 h^2 + \|{f}_m-\tilde{f}_m\|_2^2) \nonumber \\
        \le& \tau^2\left(C_1 h^2 + 2 L_f^2(1+\frac{2}{h^2})\|\mathcal{P}_h{U}_m-\tilde{U}_m\|_2^2\right) \label{eq:I21_error}
    \end{align}
    with $C_1=2C_0^2L_f^2|U_m|_{H^1}^2+2C_0^2|f_m|_{H^1}^2+4L_f^2C_u^2$ depending on the smoothness of $u$ and $f$, where we obtain the first inequality in a similar manner to Eq.~\eqref{eq:I11_error}
    and the last one using Eqs.~\eqref{eq:ff}-\eqref{eq:UU}. In particular, if $f$ is independent of $\nabla u$, the error bound can be improved to
    \begin{equation}
        \label{eq:I21_error_nabla_free}
        h^d\sum_{k=1}^K I^m_{2,1}(\varphi_k^\tau) \le \tau^2\left(C_1 h^2 +2 L_f^2\|\mathcal{P}_h{U}_m-\tilde{U}_m\|_2^2\right),
    \end{equation}
    for some constant $C_1=2C_0^2L_f^2|U_m|_{H^1}^2+2C_0^2|f_m|_{H^1}^2$.

    For the reconstruction variance part of Eq.~\eqref{eq:I2}, similar to Eqs.~\eqref{eq:I12} and \eqref{eq:I12_error}, but with the only difference that the weights of born particles are $\pm 1$ rather than $\pm \tilde{Z}_m$, we have
    \begin{equation}
        \begin{aligned}
            \label{eq:I22_error}
            I^m_{2,2}(\varphi^\tau)\le \frac{1}{N(0)}\langle (\varphi^\tau)^2, |\tau \mathcal{P}_h\tilde{f}_m|\rangle, \quad h^d\sum_{k=1}^K I^m_{2,2}(\varphi_j^\tau)\le \frac{\|\tau\mathcal{P}_h\tilde{f}_m\|_{1}}{N(0) h^d}=\frac{\mathbb{E}[N_b(t_m)\big|S(t_m)]}{N(0)^2h^d}.
        \end{aligned}
    \end{equation}

    \item {\it As for $I^m_3(\varphi^\tau)=2\mathbb{E}[(\langle \varphi^\tau, {U}_{m}-X_{t_m,o}\rangle \langle\varphi^\tau, \tau {f}_m-X_{t_m,b}\rangle )\big|S(t_m)]$}. The term is the crossing error between the death and birth mechanism. The independence of the two mechanisms implies 
        \begin{equation}
        \begin{aligned}
            I^m_3 &= 2\mathbb{E}\left[\left(\langle\varphi^\tau, {U}_m-X_{t_m,o}\rangle\right)\big|S(t_m)\right] \mathbb{E}\left[\left(\langle\varphi^\tau, \tau {f}_m - X_{t_m,b}\rangle\right)\big|S(t_m)\right]\\
            &= 2\mathbb{E}\left[\left(\langle\varphi^\tau, {U}_m-X_{t_m,o}\rangle\right)\big|S(t_m)\right]\left(\langle\varphi^\tau, \tau {f}_m - \tau \mathcal{P}_h\tilde{f}_m\rangle\right), 
        \end{aligned}
    \end{equation}
    which eliminates the variance part of the birth mechanism by the conditional expectation~\eqref{eq:birth_exp}, thereby providing a higher order of $\tau$ instead of $\sqrt{\tau}$. Then using the am-gm inequality $2xy\le\tau x^2+\frac{1}{\tau}y^2$, we are able to obtain an estimation  
    \begin{equation}
      I^m_3(\varphi^\tau)\le \tau I^m_1(\varphi^\tau) + \frac{1}{\tau}I^m_{2,1}(\varphi^\tau).
    \end{equation}
    \end{itemize}
 After finishing the central estimation of Term $I_1^m$, $I_2^m$, and $I_3^m$ in the flowchart, from Eq.~\eqref{eq:I123}, we arrive at 
     \begin{align}
       \mathbb{E}\left[e_{m+1}^2(\varphi)\big|S(t_m)\right] \le& (1+\tau) I^m_1(\varphi^\tau) + I^m_{2,2}(\varphi^\tau) + (1+\frac{1}{\tau}) I^m_{2,1}(\varphi^\tau)\nonumber \\
        \le& (1+\tau) I^m_1(\varphi^\tau) + I^m_{2,2}(\varphi^\tau) + \frac{A}{\tau}I^m_{2,1}(\varphi^\tau),  
        \label{eq:local_error}
    \end{align}
    for some constant $A\ge 1+\tau$.


    Next, we turn to cumulative error. At the first stage, we want to estimate the cumulative error between two consecutive death steps because the death mechanism is not triggered at every time step. Note that the only difference between the cases with or without the death mechanism at time $t_m$ lies in the estimation of $I^m_1$. 
    When the death mechanism occurs, we have the estimation~\eqref{eq:I11_error} for the bias part, while the variance part can be bounded by Eq.~\eqref{eq:I12_error}. When no death mechanism occurs, $I^m_1$ works as a propagation of the previous error as Eq.~\eqref{eq:I1_error_no_death}. Suppose that the death mechanism occurs $R$ times during the time interval $[0,T]$ at time steps $0<m_1<m_2<\cdots<m_R\le T/\tau$, and denote $m_0=0$. For any time step $m\in[m_r, m_{r+1})$ with $r=0,1,\dots,R$, the local error inequality~\eqref{eq:local_error} provides a recursive relation between time step $m$ and $m+1$, i.e.,
    \begin{equation}
      \label{eq:iter}
        \mathbb{E}\left[e_{m+1}^2(\varphi)\big|S(t_m)\right] \le (1+\tau) e_m^2(\varphi^\tau) + I^m_{2,2}(\varphi^\tau) + \frac{A}{\tau}I^m_{2,1}(\varphi^\tau).
    \end{equation}
    Since we have the Markovian property that $S(t_{m+1})$ depends only on $S(t_m)$, the conditional expectation on the left side of the inequality above can be written into
    \begin{equation*}
      \mathbb{E}\left[e_{m+1}^2(\varphi)\big|S(t_m)\right] = \mathbb{E}\left[e_{m+1}^2(\varphi)\big|S(t_m), S(t_{m-1}),\dots, S(t_{m_r})\right].  
    \end{equation*}
    Combining this with the tower property of conditional expectation, we have
    \begin{equation*}
        \mathbb{E}\left[e_{m+1}^2(\varphi)\big|S(t_{m_r})\right] = \mathbb{E}\left[\mathbb{E}\left[e_{m+1}^2(\varphi)\big|S(t_m), S(t_{m-1}),\dots, S(t_{m_r})\right]\big|S(t_{m_r})\right].
    \end{equation*} 
    Therefore, we can take conditional expectation on $S(t_{m_r})$ for inequality~\eqref{eq:iter}, get
    \begin{equation*}
      \mathbb{E}\left[e_{m+1}^2(\varphi)\big|S(t_{m_r})\right] \le (1+\tau) \mathbb{E}\left[e_m^2(\varphi^\tau)\big|S(t_{m_r})\right] + \mathbb{E}\left[I^m_{2,2}(\varphi^\tau)\big|S(t_{m_r})\right] + \frac{A}{\tau} \mathbb{E}\left[I^m_{2,1}(\varphi^\tau)\big|S(t_{m_r})\right], 
    \end{equation*}
    and then iterate the inequality above from time step $m_r$ to $m-1$, gaining for $m=m_r+1,\dots,m_{r+1}$, 
    \begin{equation}
        \begin{aligned}
            \mathbb{E}\left[e_{m}^2(\varphi)\big|S(t_{m_r})\right] \le& (1+\tau)^{m-m_r} I^{m_r}_1(\varphi^{(m-m_r)\tau}) + \sum_{j=m_r}^{m-1} (1+\tau)^{m-1-j} \mathbb{E}\left[I^j_{2,2}(\varphi^{(m-j)\tau})\big|S(t_{m_r})\right]\\
            &+ \frac{A}{\tau} \sum_{j=m_r}^{m-1} (1+\tau)^{m-1-j} \mathbb{E}\left[I^j_{2,1}(\varphi^{(m-j)\tau})\big|S(t_{m_r})\right].
        \end{aligned}
    \end{equation}
    Take $\varphi=\varphi_k$, sum over $k=1,2,\dots,K$, and employ the estimations~\eqref{eq:I21_error} and \eqref{eq:I22_error}, as well as the identity~\eqref{eq:L2_error}, to get
    \begin{equation}
        \begin{aligned}
            \frac{ \mathbb{E}\left[\|\mathcal{P}_h{U}_m - \tilde{U}_{m}\|_{2}^2\big|S(t_{m_r})\right]}{(1+\tau)^{m-m_r}}\le&  h^d\sum_{k=1}^K I^{m_r}_1(\varphi_k^{(m-m_r)\tau})+\frac{\mathbb{E}[\sum_{j=m_r}^{m-1}N_b(t_j)\big|S(t_{m_r})]}{N(0)^2h^d}\\
            &+AC_1h^2(1-(1+\tau)^{-(m-m_r)})\\
            &+2AL_f^2(1+\frac{2}{h^2})\tau\sum_{j=m_r}^{m-1}\frac{\mathbb{E}\left[\|\mathcal{P}_h{U}_j - \tilde{U}_j\|_2^2\big|S(t_{m_r})\right]}{(1+\tau)^{j-m_r}}.
        \end{aligned}
    \end{equation}
    Substituting 
    \begin{equation*}
        \begin{aligned}
            a_m&=\frac{\mathbb{E}[\|\mathcal{P}_h{U}_m-\tilde{U}_m\|^2\big|S(t_{m_r})]}{(1+\tau)^{m-m_r}}, \\
            c_m&=h^d\sum_{k=1}^K I^{m_r}_1(\varphi_k^{(m-m_r)\tau})+\frac{\mathbb{E}[\sum_{j=m_r}^{m-1}N_b(t_j)\big|S(t_{m_r})]}{N(0)^2h^d}+AC_1h^2(1-(1+\tau)^{-(m-m_r)}),
        \end{aligned}
    \end{equation*}
     into the discrete Gronwall inequality in Lemma \ref{lem:gronwall} and using Eqs.~\eqref{eq:L2_error}, \eqref{eq:I11_error} and \eqref{eq:I12_error}, we have for $r=1,2,\dots,R$,
    \begin{equation}
        \label{eq:between_death}
        \begin{aligned}
            &\mathbb{E}\left[\|{U}_m - \tilde{U}_{m}\|_2^2\big| S(t_{m_r})\right] \\
            \le&
            \mathbb{E}\left[\|\mathcal{P}_h{U}_m - \tilde{U}_{m}\|_{2}^2\big|S(t_{m_r})\right]+C_0|{U}_m|_{H^1}^2 h^2\\
            \le&\me^{(B+1)(t_m - t_{m_r})}\left(\|{U}_{m_r}-\tilde{U}_{m_r}\|_2^2+\frac{\tilde{Z}_{m_r}^2}{N(0)h^d}+ \frac{\mathbb{E}[\sum_{j=m_r}^{m-1}N_b(t_j)\big|S(t_{m_r})]}{N(0)^2h^d}\right)\\
            &+\me^{B(t_m-t_{m_r})}(\me^{t_m-t_{m_r}}-1)AC_1h^2 +C_0|{U}_m|_{H^1}^2 h^2,
        \end{aligned}
    \end{equation}
    and finish the local cumulative error estimation as Step 3 in Figure~\ref{fig:flowchartproof}, where $B=2AL_f^2(1+\frac{2}{h^2})$.
    
    Finally, we let $m=m_{r+1}$, iterate the inequality~\eqref{eq:between_death} globally from $r=0$ to $r=R-1$ as Step 4 in the flowchart and obtain
    \begin{equation}
        \begin{aligned}
            \mathbb{E}\left[\|{U}_m - \tilde{U}_{m}\|_{2}^2\right] \le& \me^{(B+1) T}\left(\frac{\mathbb{E}[\sum_{j=0}^{m-1}N_b(t_j)]}{N(0)^2h^d}+AC_1h^2\right.\\
            &+\left. (\mathbb{E}[R]+1)\max_{0\le k\le m}\left(\frac{\mathbb{E}[\tilde{Z}_k^2]}{N(0)h^d}+C_0|U_k|_{H^1}^2h^2\right)\right).
        \end{aligned}
    \end{equation}
    Since the death mechanism occurs when the particle number exceeds a prescribed threshold $n_A N(0)$, we have $(n_A-1)N(0)\le \sum_{j=m_r}^{m_{r+1}-1}N_b(t_j)$ and thus  the estimation of $\mathbb{E}[R]$~\eqref{eq:R} holds. 
\end{proof}

Since we can regard SPM as a special case of SPM-birth-death with death at every time step and without any birth, Theorem~\ref{thm:error_SPM} can be directly obtained following the similar procedure in the proof of Theorem~\ref{thm:error_SPM_birth_death}.

\section{Numerical experiments}\label{par bir numer}
This section conducts several numerical experiments and makes a comparison between SPM-birth-death and SPM. The two numerical examples below are adopted from \cite{Lei2025}. The particle growth threshold parameters in Algorithm~\ref{birdie alg} is set as $n_A = 3$ unless otherwise specified.

\subsection{The 1-D benchmark}
We first consider the following 1-D nonlinear model with a Cauchy data to benchmark the two methods,
\begin{equation}\label{1d non}
	\begin{aligned}
		\frac{\partial u}{\partial t} &= { \nabla u} + \Delta u + u - u^3, \\
		u(x,0) &= u_0(x) = \exp(-x^2)(1+x^4).
	\end{aligned}
\end{equation}
The relative $L^2$ error
\begin{equation}\label{rel L2}
	\mathcal{E}_2[u](t)=\frac{\left\|u^{\mathrm{num}}(x, t)-u^{\mathrm{ref}}(x, t)\right\|_2}{\left\|u^{\mathrm{ref}}(x, t)\right\|_2}
\end{equation}
is employed to measure the accuracy. The reference solutions $u^{\mathrm{ref}}$ to Eq.~\eqref{rel L2} are produced by a deterministic solver which adopts the second-order operator splitting \cite{hansen2012second}.
The 1-D benchmark in this section are implemented via Matlab with platform parameters: Intel (R) Core (TM) i7-6700 CPU (3.40 GHz), 16.0 GB of memory. We set the final time $T = 10$, the time step $\tau = 0.1$ and the side length of hybercube $h = 0.1$ in Algorithm~\ref{birdie alg} unless otherwise indicated.

\begin{figure}[hbtp]
		\centering
		\includegraphics[width=.49\textwidth]{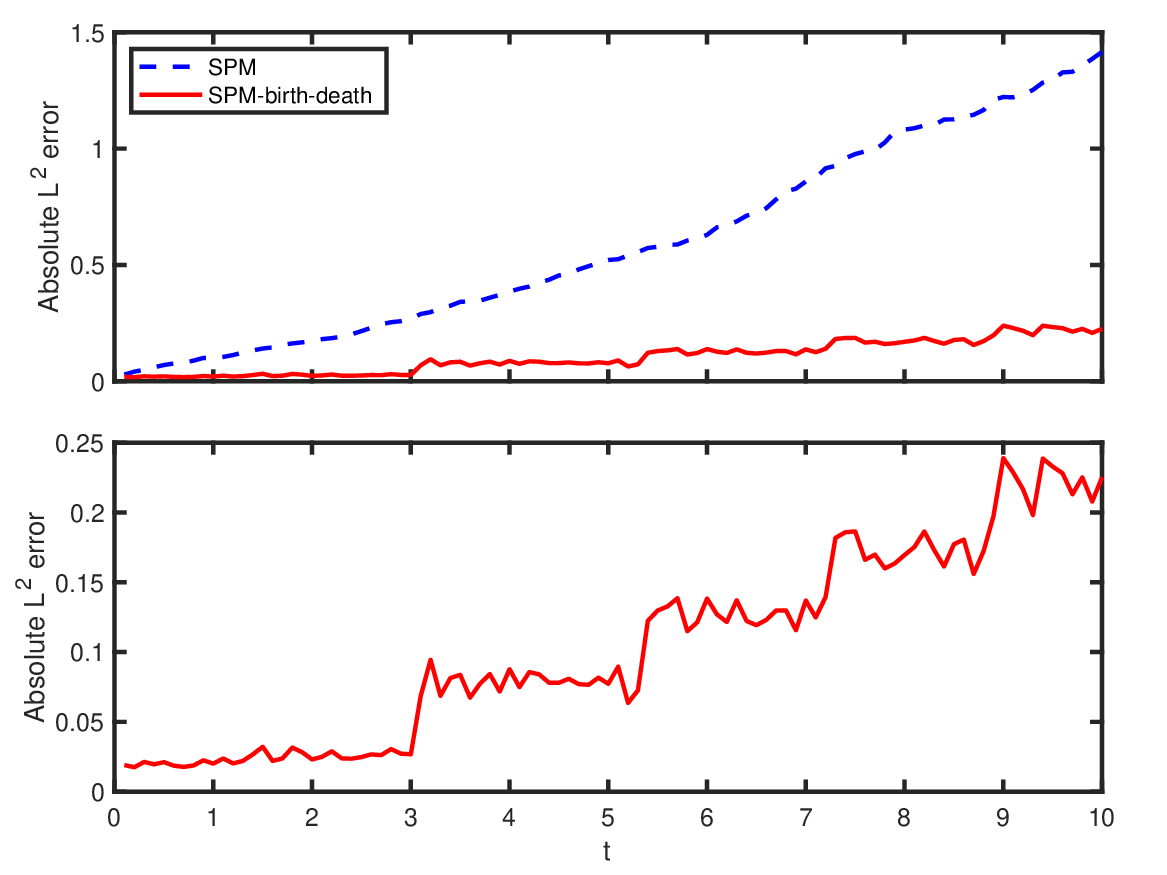}
	\caption{1-D benchmark: Results with parameters $h = 0.2$ and $\tau = 0.01$, $N=1\times 10^6$ for SPM, $N(0) = 1\times 10^5$, $n_A=9$ for SPM-birth-death. (Top) The absolute $L^2$ error of SPM and SPM-birth-death($R=4$) over time $t\in [0,10]$. The maximal particle number and memory usage of SPM-birth-death are guaranteed to be less than those of SPM since $n_AN(0)< N$. The comparison highlights the superior performance of SPM-birth-death. (Bottom) Zoomed-in view of SPM-birth-death error. The error jumps when death mechanism occurs and oscillates more thereafter, since the resampling procedure introduces the spatial bias and statistical variance. }%
	\label{fig:na9_L4_abserr}%
\end{figure}
\begin{figure}[hbtp]
	\centering
	\begin{subfigure}[b]{0.49\textwidth}
		\centering
		\includegraphics[width=\textwidth]{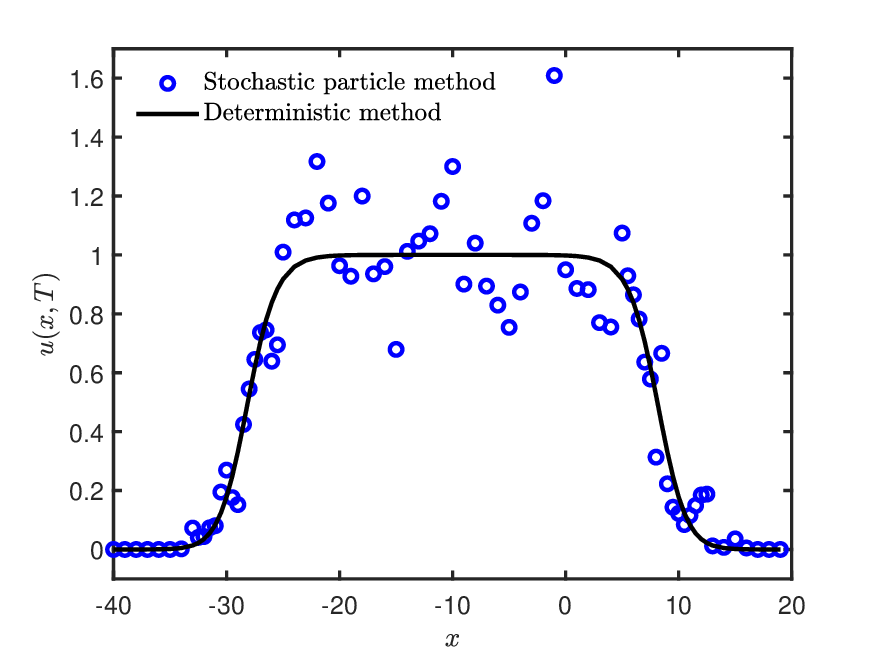}
		\caption{$N(0) = 1\times 10^4$, SPM-birth-death.}
	\end{subfigure}
	\begin{subfigure}[b]{0.49\textwidth}
		\centering
		\includegraphics[width=\textwidth]{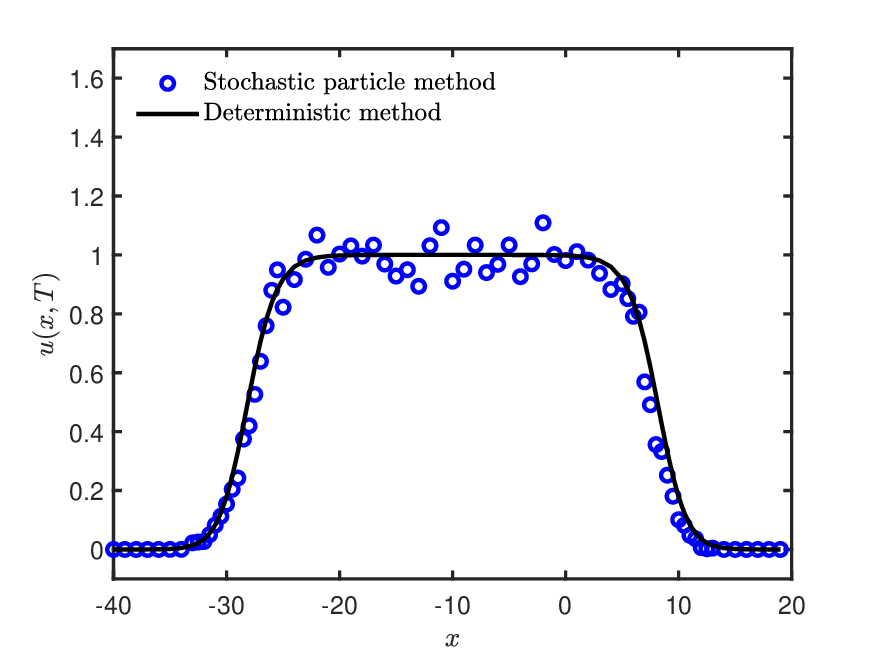}
		\caption{$N(0) = 1\times 10^5$, SPM-birth-death.}
	\end{subfigure}
	\hfill
	\begin{subfigure}[b]{0.49\textwidth}
		\centering
		\includegraphics[width=\textwidth]{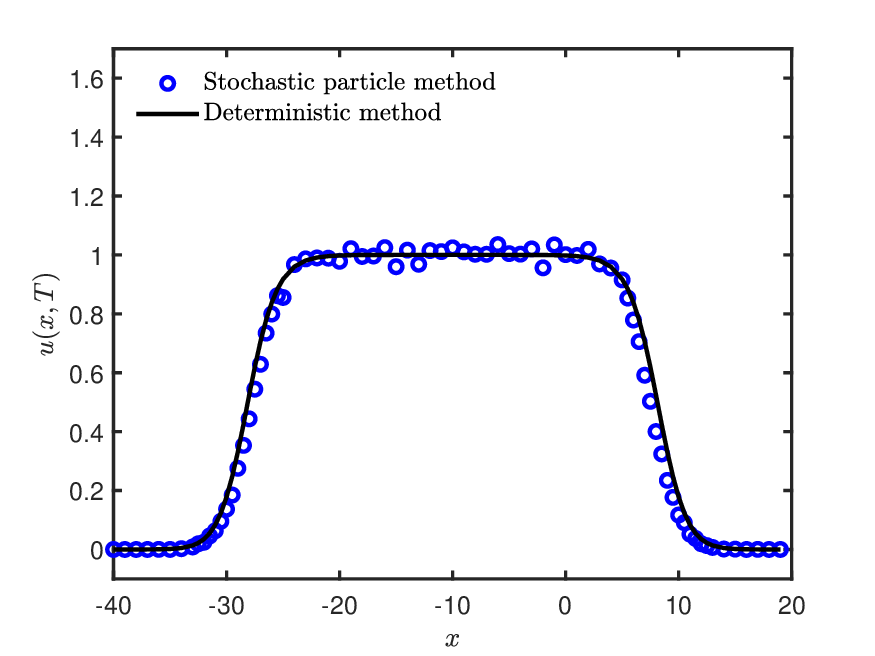}
		\caption{$N(0) = 1\times 10^6$, SPM-birth-death.}
	\end{subfigure}
	\hfill
	\begin{subfigure}[b]{0.49\textwidth}
		\centering
		\includegraphics[width=\textwidth]{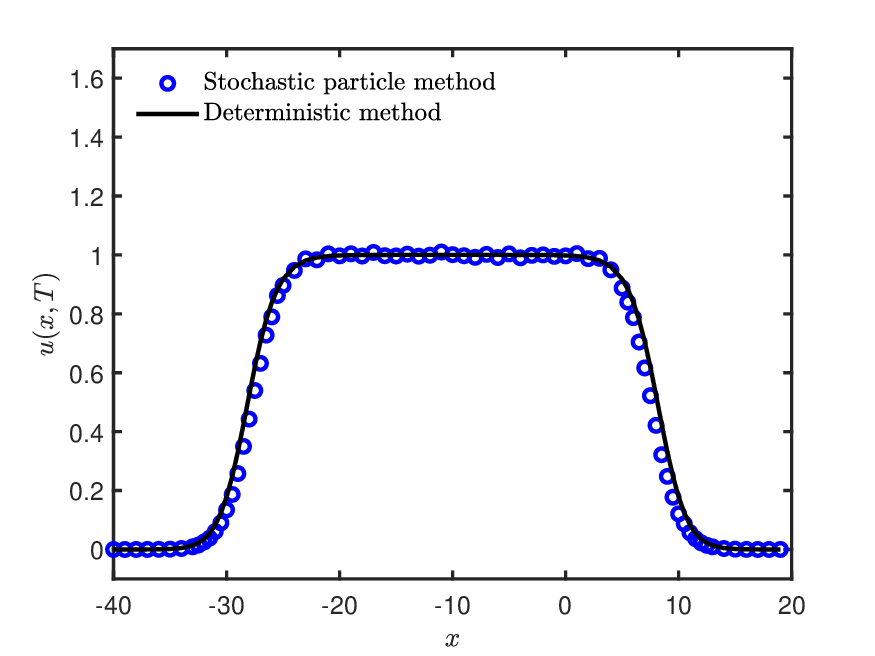}
		\caption{$N(0) = 1\times 10^7$, SPM-birth-death.}
	\end{subfigure}
	\caption{1-D benchmark: The numerical effects of SPM-birth-death with different particle numbers show that its variance decreases as the initial sampling number $N(0)$ increases.}
	\label{Fig D}
\end{figure}
\begin{figure}[hbtp]
	\centering
	\includegraphics[width=0.49\linewidth]{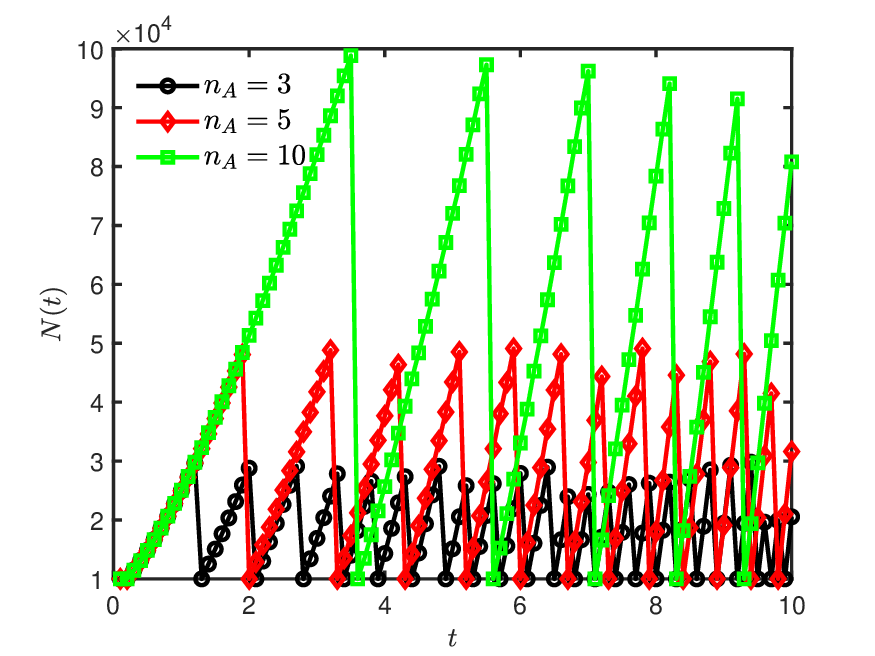}
	\caption{1-D benchmark: The variation of particle number over time for SPM-birth-death. When the particle number surpasses 
	$n_A\times N(0)$, the particle annihilation mechanism is activated, resulting in a reduction of the particle number.}
	\label{fig Dparnum}
\end{figure}
\begin{figure}[hbtp]
	\centering
	\includegraphics[width=0.49\linewidth]{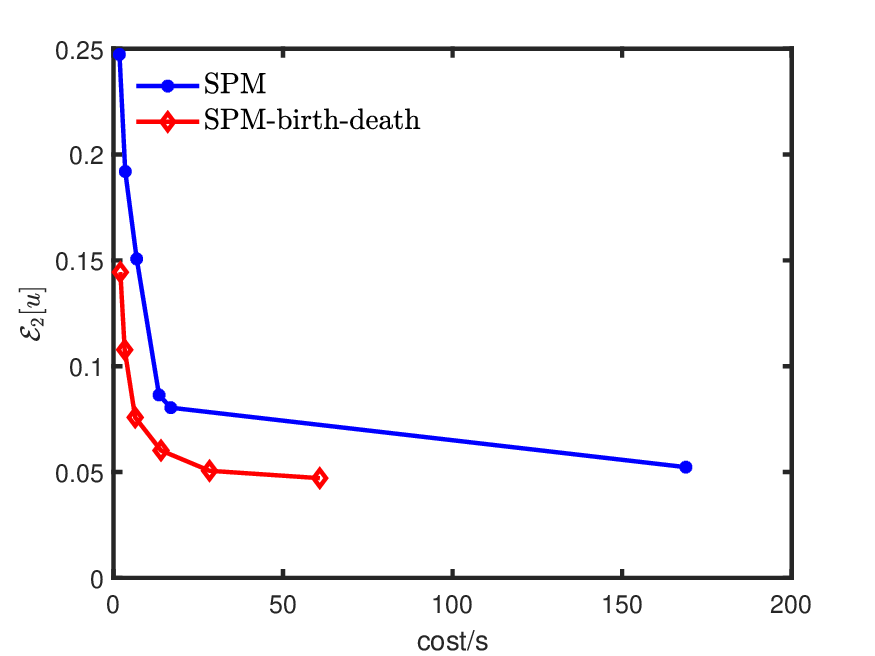}
	\caption{1-D benchmark: 
		The correlation between relative $L^2$ error $\mathcal{E}_2[u]$ and computational cost for SPM and SPM-birth-death. SPM-birth-death demonstrates superior computational efficiency: achieving lower errors at equivalent computational times (CPU time) and requiring shorter times to reach equivalent error levels. We set in SPM:  $N = 1\times 10^5, 2\times 10^5, 4\times 10^5, 8\times 10^5, 1\times 10^6, 1\times 10^7$ with the corresponding running time respectively being 1.7502, 3.4546, 6.8315, 13.42, 16.8924, 168.81 seconds, 
		and in SPM-birth-death: $N(0) = 2\times 10^4, 4\times 10^4, 8\times 10^4, 2\times 10^5, 4\times 10^5, 8\times 10^5$ with the corresponding running time respectively being 2.0102, 3.2941, 6.3855, 13.9954, 28.31, 60.82 seconds. The maximal particle number and memory usage of SPM-birth-death are guaranteed to be less than those of SPM since $n_AN(0)\le N$ for $n_A=3$. Here sample sizes were strategically selected to ensure comparable runtimes, enabling direct comparison within the same figure. }
	\label{fig BDeffici}
\end{figure}

Table~\ref{tableaccuracy} demonstrates the convergence order of SPM-birth-death. As Theorem~\ref{thm:error_SPM_birth_death} expected, the method achieves half-order accuracy with respect to initial sample size $N(0)$. Due to its Lawson-Euler scheme and piecewise constant reconstruction, it exhibits first-order convergence in both time step $\tau$ and spatial resolution $h$. Table~\ref{tab:na_converge} illustrates SPM-birth-death's error scaling as $\mathcal{O}(\sqrt{\frac{1}{n_A-1}})$ when $n_A$ approaches 1, as shown in Theorem~\ref{thm:error_SPM_birth_death}. Figure~\ref{fig:na9_L4_abserr} depicts absolute $L^2$-error over time for large $n_A$ values with less frequent annihilation, where we can see the jumps at annihilation points and growing oscillations thereafter, consistent with Theorem~\ref{thm:error_SPM_birth_death}'s analysis. 
Figure~\ref{Fig D} reveals systematic accuracy improvement under SPM-birth-death, with solution variance decreasing as $N(0)$ increases. Figure~\ref{fig Dparnum} illustrates the particle annihilation mechanism for different $n_A$ parameters: Particles are annihilated when the count exceeds $n_A \times N(0)$, reducing computational costs. Figure~\ref{fig BDeffici} confirms SPM-birth-death's superior computational efficiency versus the standard SPM \cite{Lei2025}, achieving lower errors at equivalent computational costs.


\begin{table}[htbp]
	\centering
	\caption{1-D benchmark: SPM-birth-death has first-order accuracy in both time and space plus half-order accuracy in the initial sample size.}
	\begin{tabular}{ccccccc}
		\hline 
		$N(0)$ & $h$   & $\tau$  &    $\mathcal{E}_2[u]$(10) & order  \\ 
		\midrule
		$1\times 10^4$ &  &        & 0.5783 & -    &   \\ 
		$4\times 10^4$ & 0.01 & 0.01       & 0.2988& $0.48$   \\ 
		$1.6\times 10^5$ & &     & 0.1504 & $0.49$   \\ \midrule
		& & 0.25 & 0.1023 & - \\
		$1\times 10^7$ & 0.01 & $0.2$ &  0.0809 & 1.05   \\ 
		&  & $0.1$   & 0.0414 & $0.99$    \\ \midrule
		& $0.2$ &       & 0.0106 & -    \\ 
		
		$1\times 10^7$ & $0.15$ & 0.01      & 0.0076 & 1.15       \\ 
		& $0.1$ &        & 0.0051 & $1.05$    \\
		\hline 
	\end{tabular}
	\label{tableaccuracy}
\end{table}

\begin{table}[htbp]
  \centering
  \caption{1-D benchmark: SPM-birth-death has the error of $\mathcal{O}(\sqrt{\frac{1}{n_A-1}})$ when $n_A$ is close to 1, for fixed $N(0) = 1\times 10^6$, $h = 0.2$ and $\tau = 0.01$. Here $R$ denotes the number of death mechanisms occurred during $t\in[0,10]$.}
    \begin{tabular}{cccccc}
	\hline
    $n_A$    & $\frac{1}{n_A-1}$ & $\mathcal{E}_2[u]$(10) & R & order \\
	\midrule
	1.1   & 10    & 0.0560 & 302 & - \\
	1.2   & 5     & 0.0398 & 163 & 0.49 \\
    1.5   & 2     & 0.0242 & 67  & 0.55 \\
	\hline
    \end{tabular}%
  \label{tab:na_converge}%
\end{table}%

%

\subsection{The d-D Allen-Cahn equation}
\par
\begin{figure}[hbtp]
	\centering
	\begin{subfigure}[b]{0.3\textwidth}
		\centering
		\includegraphics[width=\textwidth]{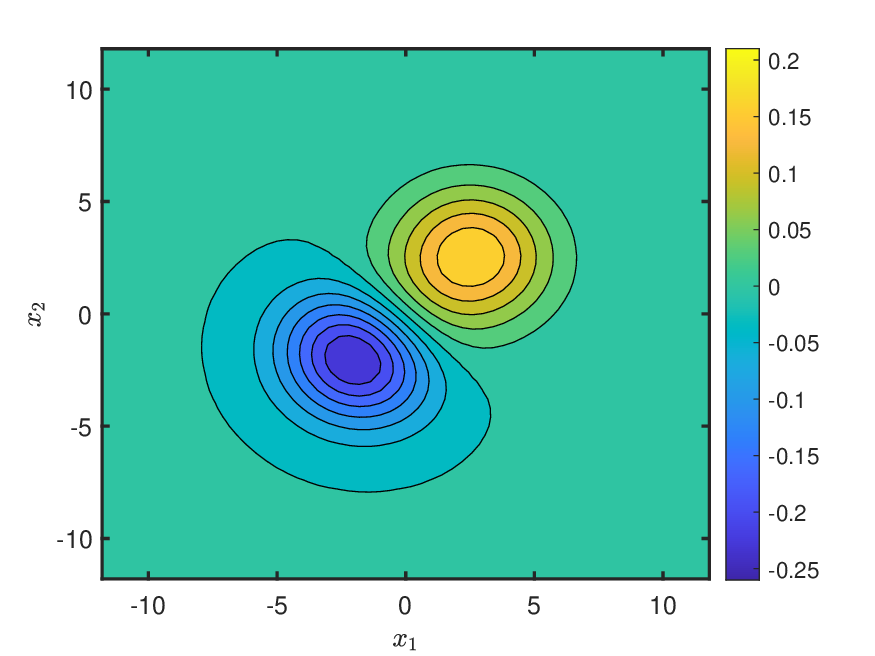}
		\caption{Reference solution.}
	\end{subfigure}
	\hfill
	\begin{subfigure}[b]{0.3\textwidth}
		\centering
		\includegraphics[width=\textwidth]{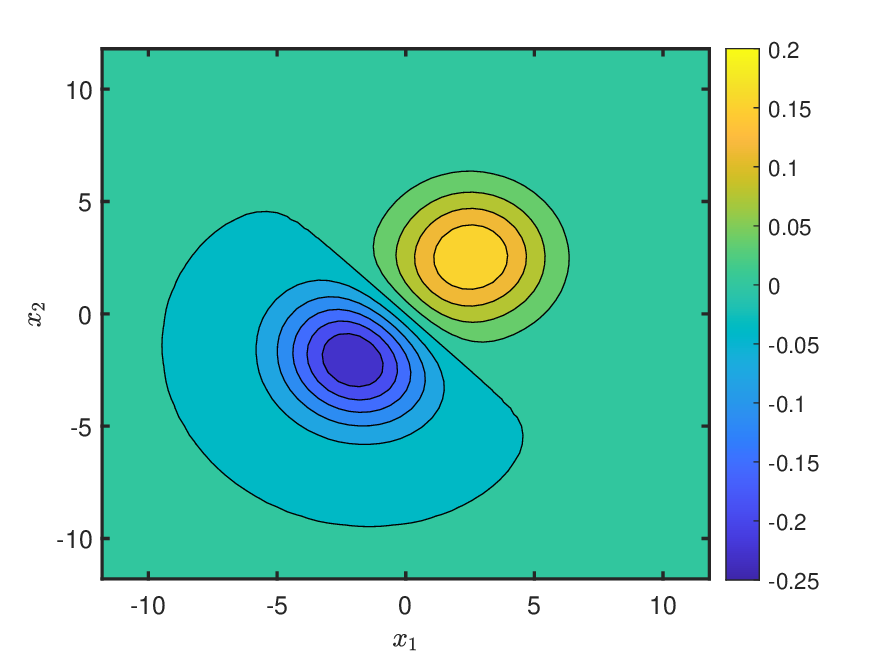}
		\caption{SPM-birth-death, $N(0)=1\times 10^8$.}
	\end{subfigure}
	\hfill
	\begin{subfigure}[b]{0.3\textwidth}
		\centering
		\includegraphics[width=\textwidth]{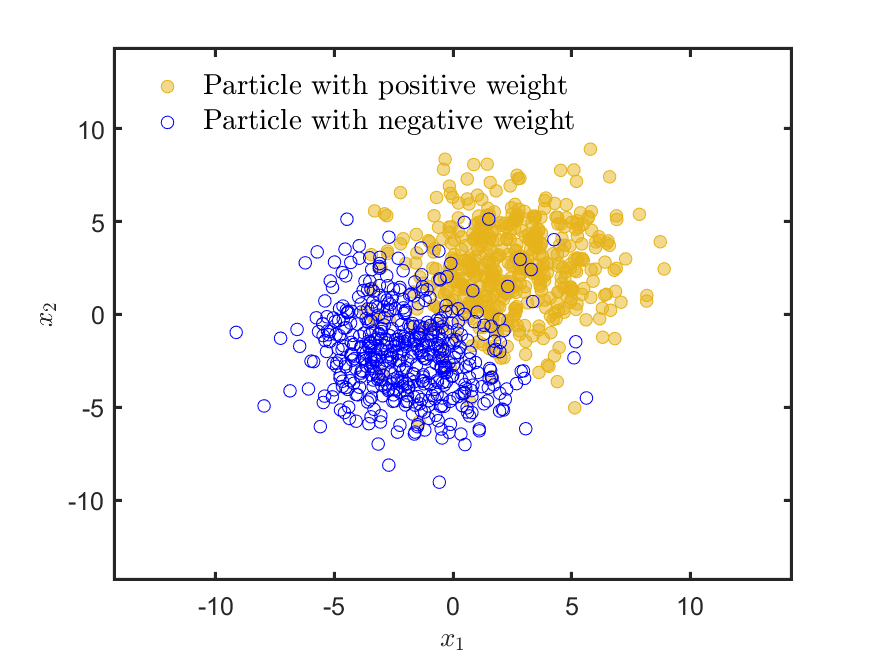}
		\caption{Particle distribution. \label{fig Dallencahnparticlex}}
	\end{subfigure}
	\caption{The 6-D Allen-Cahn equation: Filled contour plots of $M(x_1, x_2, 2)$ defined in Eq.~\eqref{P&M} for (a) the reference solution given in
		Eq.~\eqref{allencahnexact} and (b) the numerical solution produced by SPM-birth-death with $N(0)=1\times 10^8$. We randomly choose $10^3$ particles from all samples at the final instant $T=2$ and project their locations onto the $x_1x_2$-plane in (c) to show that the particles in SPM-birth-death exhibit highly adaptive characteristic. The particles with positive and negative weights concentrate in the areas where the solution takes positive and negative values, respectively.}
	\label{fig Dallencahn}
\end{figure}

\begin{figure}[hbtp]
	\centering
	\begin{subfigure}[b]{0.3\textwidth}
		\centering
		\includegraphics[width=\textwidth]{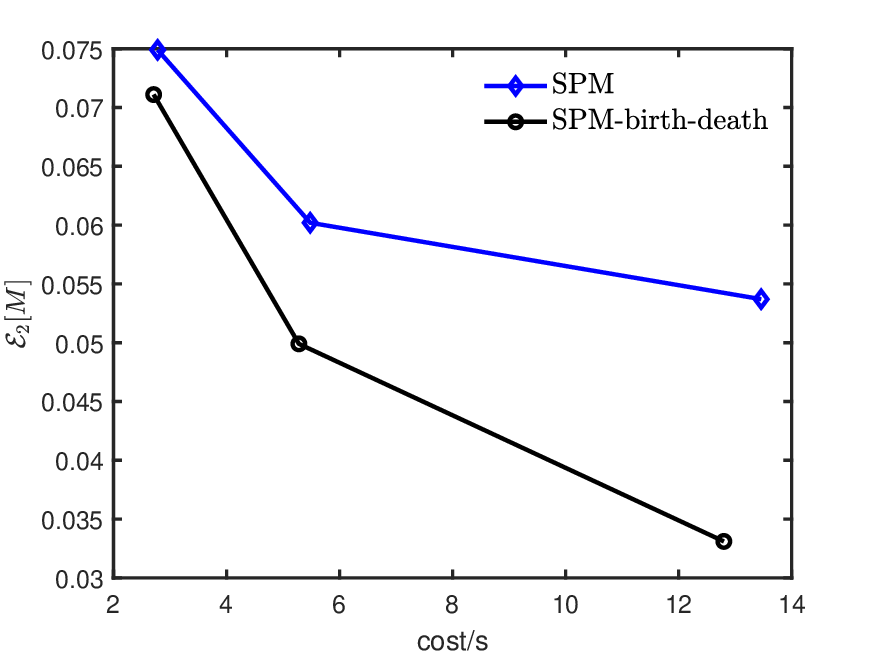}
		\caption{2-D.}
	\end{subfigure}
	\hfill
	\begin{subfigure}[b]{0.3\textwidth}
		\centering
		\includegraphics[width=\textwidth]{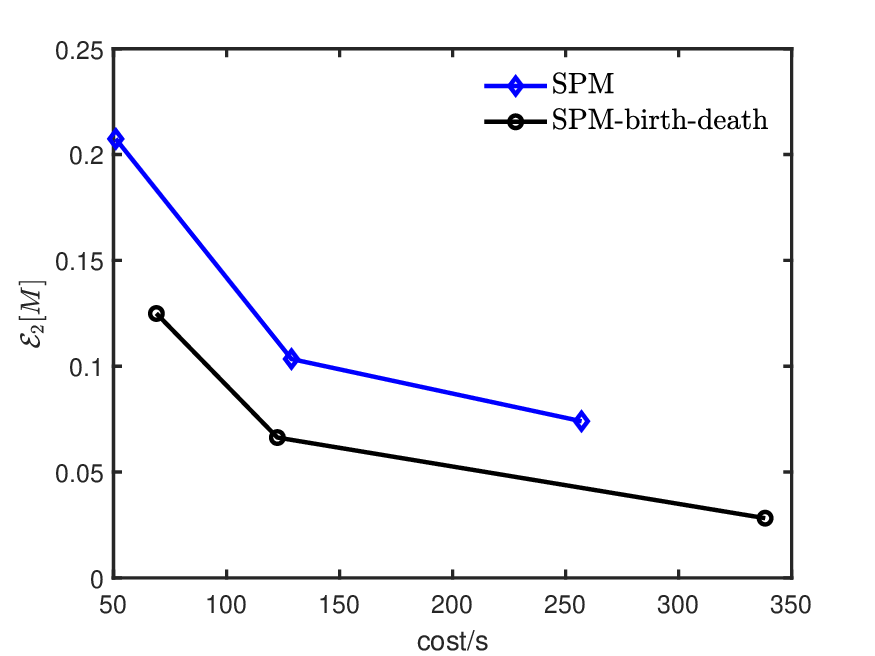}
		\caption{4-D.}
	\end{subfigure}
	\hfill
	\begin{subfigure}[b]{0.3\textwidth}
		\centering
		\includegraphics[width=\textwidth]{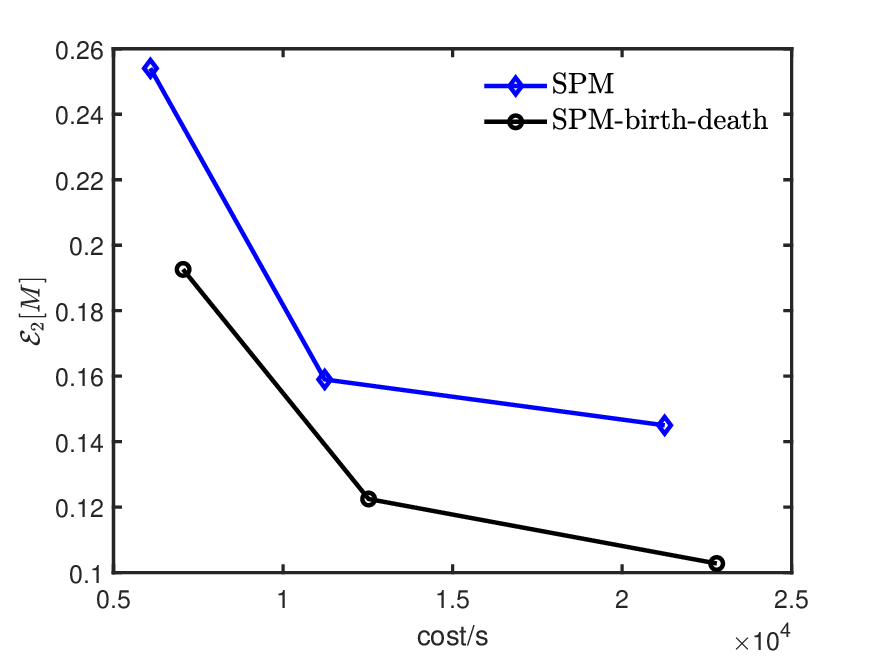}
		\caption{6-D.}
	\end{subfigure}
	\caption{The d-D Allen-Cahn equation: The relationship between relative $L^2$ error $\mathcal{E}_2[M]$ and computational time cost for SPM and SPM-birth-death in different dimensions. Compared with SPM, SPM-birth-death can achieve smaller errors at the same computational cost in various dimensions. By selecting different sample sizes in the two methods to ensure comparable running times, we can plot their results in the same figure. Here we set the side length of hybercube $h = 0.3$ in 2-D and 4-D experiments.}
	\label{fig 246dallencahn}
\end{figure}

 Next we consider the high-dimensional numerical experiments with the $d$-D Allen-Cahn equation, 
\begin{equation}\label{ac}
	\frac{\partial}{\partial t}u(\vec{x},t) = c\Delta u(\vec{x},t) +  u(\vec{x},t)-{ u^3(\vec{x},t)} + {r(\vec{x},t)},\quad \bx \in \mathbb{R}^d,
\end{equation}
which allows the following analytical solution
\begin{equation}\label{allencahnexact}
	u^{\text{ref}}(\vec{x},t) = \frac{x_1+x_2}{\left(\pi (1+4ct)\right)^{d/2}}\left(\exp\left(-\frac{||\vec{x}-\vec{p}_1||_2^2}{1+4ct}\right) +2.0 \exp\left(-\frac{||\vec{x}-\vec{p}_2||_2^2}{1+4ct}\right)\right),
\end{equation}
where $\vec{p}_1 = \left(2,2,0,\dots,0\right)$, $\vec{p}_2 = \left(-1,-1,0,\dots,0\right)$ and $c = 1$ is the diffusion coefficient. The initial data $u_0(\bx) = u^{\text{ref}}(\vec{x},0)$ and let $r(\vec{x},t)$ be the remaining term after substituting the analytical solution Eq.~\eqref{allencahnexact} into Eq.~\eqref{ac}. We set the final time $T = 2$, the time step $\tau = 0.1$ and the side length of hybercube $h = 0.4$ in Algorithm~\ref{birdie alg} unless otherwise specified. To visualize high-dimensional solutions, we adopt the following 2-D projection,
\begin{equation}\label{P&M}
	 M(x_1,x_2,t) = \int_{\mathbb{R}^{d-2}} u(\bx,t) \D x_3 \dots \D x_d,
\end{equation}
and use the relative $L^2$ error $\mathcal{E}[M](t)$ to measure the accuracy. The 2-D projection is approximated by a uniform partition (take test function $\varphi(\bx) = \mone_{Q^{\mu}\times Q^{\nu}}$),
\begin{equation}\label{lowdim proj}
		M(x_1,x_2, t)  \approx \sum_{\mu=1}^{(r_1-l_1)/h} \sum_{\nu=1}^{(r_2-l_2)/h} \left(\frac{1}{N(t) |Q^{\mu}| |Q^{\nu}|} \sum_{i=1}^{N(t)} w_i(t) \mone_{Q^{\mu}\times Q^{\nu}} (\bx_i(t))\right) \mone_{Q^{\mu}\times Q^{\nu}} (x_1, x_2),
\end{equation}
where $Q^{\mu} = \left[l_1+(\mu-1)h, l_1+\mu h\right]$, $Q^{\nu} = \left[l_2+(\nu-1)h, l_2+\nu h\right]$. The simulations in this section via our C++ implementations run on the High-Performance Computing Platform of Peking University: 2*Intel Xeon E5-2697A-v4 (2.60 GHz, 40 MB Cache, 9.6GT/s QPI Speed, 16 Cores, 32 Threads) with 256 GB Memory $\times$ 16. 
Table~\ref{table:2dallen_h} demonstrates that SPM-birth-death achieves first-order accuracy in space for the 2-D Allen-Cahn equation, as expected in Theorem~\ref{thm:error_SPM_birth_death}. Figure~\ref{fig Dallencahn} plots $M(x_1,x_2, 2)$ defined in Eq.~\eqref{P&M} produced by SPM-birth-death with $N(0) = 1\times 10^8$ against the reference solution given in Eq.~\eqref{allencahnexact}. The agreement between them is evident. The characteristic of moving with the solution for particle locations in SPM-birth-death is clearly demonstrated in the last plot of Figure~\ref{fig Dallencahn}. Figure~\ref{fig 246dallencahn} shows a comparison of the computational efficiency between SPM and SPM-birth-death in the high-dimensional example. Similar to the 1-D benchmark tests, SPM-birth-death is more efficient than SPM owing to the birth-death mechanism.

\begin{table}[htbp]
	\centering
	\caption{The 2-D Allen-Cahn equation: SPM-birth-death has first-order accuracy in space.}
	\begin{tabular}{ccccccc}
		\hline 
		$N(0)$ & $h$   & $\tau$  &    $\mathcal{E}_2[u]$(2) & order  \\ 
		\midrule
		& $0.4$ &       & 0.0237 & -    \\ 
		
		$1\times 10^7$ & $0.35$ & 0.01      & 0.0203 & 1.15       \\ 
		& $0.3$ &        & 0.0171 & $1.12$    \\
		\hline 
	\end{tabular}
	\label{table:2dallen_h}
\end{table}

%

\section{Conclusion and discussion}

The primary contribution of this work is the proposal of SPM-birth-death—a method that inherits the core framework and design philosophy of SPM [34] but incorporates an active particle birth-death mechanism to enhance computational efficiency. A rigorous error estimation reveals that both methods achieve first-order accuracy in time and space, and half-order accuracy in the initial sample size with explicit variance estimates. This serves also the first theoretical justification of SPM for the existing numerical convergence study. Numerical experiments across one-dimensional and high-dimensional settings validate the convergence order and demonstrate SPM-birth-death's superior accuracy over standard SPM at equivalent computational costs. Several aspects remain open for exploration, including:
(i) alternative birth-death strategies, especially adaptively adjusting the threshold $n_A$;
(ii) more efficient methods for function reconstruction from particle clouds, and (iii) rigorous quantitative error analysis for $\mathcal{L}=\bb\cdot\nabla+c\Delta$ and more general linear operators.

\section*{Acknowledgement}
This research was supported by the National Natural Science Foundation of China (Nos.~12325112, 12288101) and the High-performance Computing Platform of Peking University. The authors are thankful to the anonymous referees for their  valuable suggestions.



\end{document}